\documentclass[11pt%
]{article}

\usepackage{amsmath}
\usepackage{amssymb}
\usepackage{amsthm}

\usepackage{graphics}
\usepackage{epigraph}
\usepackage{soul}

\usepackage{smart}
\SectioningParameter{presecnum}{section}=\noexpand\S;
\RestoreLaTeXeqno

\usepackage[all]{xy}

\usepackage{lscape}
\usepackage{rotating}

\usepackage{bbm}

\newcounter{notes}
\newenvironment{arab}{%
 \begin{list}{\arabic{notes})}{\usecounter{notes}%
  \settowidth{\labelwidth}{0.0em}
  \setlength{\labelsep}{0.0em}
  \setlength{\leftmargin}{0em}
  \setlength{\itemindent}{0.0em}
  \setlength{\itemsep}{1mm}
  \setlength{\topsep}{2mm}
  \setlength{\parsep}{3mm}
  \setlength{\partopsep}{0mm}
  }
 }
{\end{list}}

\newcommand {\supplus}{\mathop{{\supset}\llap{\raise
0.5pt\hbox{\normalfont\small+}\hskip 0.5pt}}}

\newcommand {\subplus}{\mathop{{\subset}\llap{\raise
0.5pt\hbox{\normalfont\small+}\hskip 0.5pt}}}

\newcommand {\divby}  {\lower 0.15ex \hbox{\,\vdots\,}}
\newcommand {\Cee}    {{\mathbb  C}}

\newcommand {\Kee}    {{\mathbb  K}}
\newcommand {\Nee}    {{\mathbb  N}}

\newcommand {\Zee}    {{\mathbb  Z}}

\newcommand {\fa}     {{\mathfrak{a}}}

\newcommand {\fb}     {{\mathfrak{b}}}
\newcommand {\fc}    {{\mathfrak{c}}}

\newcommand {\fe}     {{\mathfrak{e}}}

\newcommand {\fg}     {{\mathfrak{g}}}    %
\newcommand {\fgl}    {{\mathfrak{gl}}}  %
\newcommand {\fh}     {{\mathfrak{h}}}
\newcommand {\fhei}   {{\mathfrak{hei}}}
\newcommand {\fii}    {{\mathfrak{i}}}    %
\newcommand {\fn}     {{\mathfrak{n}}}

\newcommand {\fo}     {{\mathfrak{o}}}
\newcommand {\fosp}   {{\mathfrak{osp}}}

\newcommand {\fr}     {{\mathfrak{r}}}

\newcommand {\fs}     {{\mathfrak{s}}}

\newcommand {\fsl}    {{\mathfrak{sl}}}

\renewcommand {\cal} {\mathcal}

\newcommand {\cO}     {{\cal O}}

\def \opname#1#2%
  {\expandafter\newcommand \csname #1\endcsname {{\mathop{\mathrm{#2}}\nolimits }}}

\newcommand{\rmname}[1]
  {\expandafter\newcommand \csname #1\endcsname {{\operatorname{#1}}}}

\newcommand{\rmnameii}[2]
  {\expandafter\newcommand \csname #1\endcsname {{\operatorname{#2}}}}

\rmname{act} \rmname{Ad} \rmname{Add} \rmname{ad} \rmname{Alt}
\rmname{alt} \rmname{Ann} \rmname{antidiag} \rmname{Ber}
\rmname{ber} \rmname{Bil} \rmname{Br} \rmname{card} \rmname{ch}
\rmname{Char} \rmname{cem} \rmname{cj} \rmname{Cliff}
\rmname{cntr} \rmname{codim} \rmname{Coind} \rmname{const}
\rmname{col} \rmname{cork} \rmname{cpr} \rmname{diag}
 \rmnameii{Div}{div} \rmname{Def} \rmname{Deg}
\rmname{Der} \rmname{Diff} \rmname{Dim} \rmname{End} \rmname{Even}
\rmname{Ext} \rmname{gr} \rmname{Hom} \rmname{HT}
\rmnameii{Ht}{ht} \rmname{hwt} \rmname{Id} \rmname{id}
\rmname{ind} \rmname{Ind} \rmname{Inf} \rmname{irr} \rmname{Le}
\rmname{Lie} \rmname{lwt} \rmname{mult} \rmname{Mat} \rmname{Mor}
\rmname{nm} \rmname{Ob} \rmname{Odd} \rmname{Osc} \rmname{per}
\rmname{Pic} \rmname{pr} \rmname{pro} \rmname{Prime} \rmname{Proj}
\rmname{prt} \rmname{pt} \rmname{Q} \rmname{qet} \rmname{qtr}
\rmname{rd} \rmname{rk} \rmname{row} \rmname{Res} \rmname{salt}
\rmname{Sch} \rmname{SBr} \rmname{sdim}\rmname{scalar}
\rmname{Ser} \rmname{sign} \rmname{Smbl} \rmname{spin}
\rmname{ssym} \rmname{str} 
\rmname{sgn} \rmname{sq}
\rmname{symm} \rmname{supp} \rmname{Supp} \rmname{St}
\rmname{Spec} \rmname{Spm} \rmname{tr} \rmname{vpt} \rmname{Vect}
\rmname{weyl} \rmname{Weyl} \rmname{Witt}

\opname{vvol}  {{v\hspace{-0.1ex}o\hspace{-0.02ex}l\/}}
\opname{pnt}  {\text{\normalfont pt}} \opname{Span} {{Span}}
\opname{slim} {\overline{\lim}} \opname{Vol}
{{V\hspace{-0.55ex}o\hspace{-0.02ex}l\/}} \opname{Par}
{{P\hspace{-0.3ex}a\hspace{-0.05ex}r\/}}

\newcommand {\ev} {{\bar0}}
\newcommand {\od} {{\bar1}}

\newcommand {\bcdot}   {\mathbin{\hbox{\raise.4ex\hbox{\bf.}}}} 

\newcommand{\ffcirc}{{\fontsize{7.5}{9pt}\selectfont \text{\raisebox{-1.5ex}[0pt]{$\bigcirc$}}}}

\newcommand {\secno} {}

\makeatletter

\theoremstyle{plain}
\@@newtheorem*{Theorem}{\secno Theorem}
\@@newtheorem*{Lemma}{\secno Lemma}
\@@newtheorem*{Proposition}{\secno Proposition}
\@@newtheorem*{Corollary}{\secno Corollary}
\@@newtheorem*{Statement}{\secno Statement}
\@@newtheorem*{Problem}{\secno Problem}
\@@newtheorem*{Question}{\secno Question}
\@@newtheorem*{Conjecture}{\secno Conjecture}

\theoremstyle{definition}
\@@newtheorem*{Example}{\secno Example}
\@@newtheorem*{Examples}{\secno Examples}
\@@newtheorem*{Convention}{\secno Convention}
\@@newtheorem*{Comment}{\secno Comment}

\theoremstyle{remark}

\@@newtheorem*{Remark}{\secno Remark}
\@@newtheorem*{Remarks}{\secno Remarks}
\@@newtheorem*{Exercise}{\secno Exercise}
\@@newtheorem*{Solution}{\secno Solution}
\@@newtheorem*{Hint}{\secno Hint}

\newenvironment{lmatrix}{%
 \small\arraycolsep=4pt \left(%
  \matrix@check\lmatrix\env@matrix
}{
  \endmatrix\right)%
}

\newcommand{\?}{\nobreak\hskip.145em\nobreak\hskip\z@skip}

\renewcommand\appendix{\par
  \setcounter{section}{0}%
  \setcounter{subsection}{0}%
  \gdef\thesection{\appendixname\kern1ex\@Alph\c@section}}

\makeatother

\newcommand{\ssec}[2]{\subsection*{%
 \refstepcounter{subsection}%
 \boldmath\thesubsection.\kern1ex#2.}%
 \label{ss#1}%
}

\newcommand{\mbullet}{{\fontsize{16}{9pt}\selectfont
\text{\raisebox{-1.3pt}[0pt]{$\bullet$}}}}
\newcommand{\mcirc}{{\fontsize{16}{9pt}\selectfont
\text{\raisebox{-1.3pt}[0pt]{$\circ$}}}}
\newcommand{\motimes}{{\fontsize{10}{9pt}\selectfont
\text{\raisebox{0.3pt}[0pt]{$\otimes$}}}}

\begin{document}

\graphicspath{{./figs/}}

\title
{Cartan matrices and presentations of Cunha and Elduque
superalgebras}

\author{Sofiane Bouarroudj${}^1$, Pavel Grozman${}^2$, Dimitry Leites${}^3$
\thanks{DL is thankful to MPIMiS, Leipzig, for financial support and most
creative environment; to A.~Elduque for comments; to O.~Shirokova
for a \TeX pert help and to A.~Protopopov for a help with the
graphics; his drawing program is now available \cite{Pro}.} }

\address{${}^1$Department of Mathematics, United Arab Emirates University, Al
Ain, PO. Box: 17551; Bouarroudj.sofiane@uaeu.ac.ae\\
${}^2$Equa Simulation AB, Stockholm, Sweden; pavel@rixtele.com\\
${}^3$MPIMiS, Inselstr. 22, DE-04103 Leipzig, Germany\\
on leave from Department of Mathematics, University of Stockholm,
Roslagsv. 101, Kr\"aft\-riket hus 6, SE-104 05 Stockholm,
Sweden; mleites@math.su.se, leites@mis.mpg.de}

\keywords {Cartan prolongation, Tanaka-Shchepochkina
prolongations, nonholonomic manifold, Elduque superalgebra, Lie
superalgebra}

\subjclass{17B50, 70F25}

\maketitle

\begin{abstract} All inequivalent Cartan matrices (in other words, inequivalent
systems of simple roots) of the ten simple exceptional finite
dimensional Lie superalgebras in characteristic 3, recently
identified by Cunha and Elduque as constituents of Elduque's
superization of the Freudenthal Magic Square, are listed together
with defining relations between analogs of their Chevalley
generators.
\end{abstract}


\section{Introduction}
Recently Strade had published a monograph \cite{S} summarizing the
description of newly classified simple finite dimensional Lie
algebras over the algebraically closed fields $\Kee$ of
characteristic $p>3$, and also gave an overview of the examples
(due to Brown, Frank, Ermolaev and Skryabin) of simple finite
dimensional Lie algebras for $p=3$ with no counterparts for $p>3$.
Several researchers started afresh to work on the cases where
$p=2$ and 3, and, after a period of quietness after \cite{KL}, new
examples of simple Lie algebras with no counterparts for $p\neq 2,
3$ started to appear (\cite{J1,GL4,GG,Le1}). The $\Zee$-graded of
the above mentioned exceptional examples of simple Lie algebras
for $p=3$ (Brown, Frank, Ermolaev and Skryabin algebras) were
interpreted as algebras of vector fields preserving certain
distributions (\cite{GL4}).

Among the main points made in \cite{KL}, one was indication to
similarity between modular Lie algebras and Lie superalgebras
(even over $\Cee$). Classification of simple Lie {\it
super}algebras for $p>0$ and the study of their representations
is, however, of independent interest. A conjectural list of simple
finite dimensional Lie superalgebras over an algebraically closed
field of characteristic $p>5$, known for some time, was recently
cited in \cite{BjL}: the idea is, as in the Kostrikin--Shafarevich
conjecture, for the two classes of simple complex Lie
superalgebras $\fg$ (finite dimensional and vectorial), select a
$\Zee$-form $\fg_\Zee$ of $\fg$, if any such exists, and take
$\fg_\Kee:=\fg_\Zee\otimes _\Zee\Kee$, take a simple finite
dimensional subquotient $\fs\fii(\fg_\Kee)$ of $\fg_\Kee$;  we
call such examples {\it KSh-type algebras}. Finally, take
deformations\footnote{It is not clear, actually, if the
conventional notion of deformation can always be applied if $p>0$
(for the arguments, see \cite{LL}); to give the correct (or,
better say, universal) notion is an open problem, but let it pass
for the moment: the conventional notion is applicable to the
algebras $\fg(A)_\Kee$.} thereof if any exist.

Alberto Elduque worked lately on getting examples of simple Lie
superalgebras for $p=3, 5$ different from KSh-type. His approach
to superization of the Freudenthal Magic Square in terms of {\it
symmetric compositions algebras} led him and his Ph.D. student,
Isabel Cunha, to discovery of several new simple finite
dimensional Lie superalgebras with Cartan matrix, which,
conjecturally, are indigenous to $p=3$, cf. \cite{CE,El1,CE2}.

Interestingly, the exceptional simple Lie superalgebras
$\fa\fg(2)$ and $\fa\fb(3)$ do not appear in the superizations of
Freudenthal's magic square for $p=0$ but their analogs do appear
after their integer structure constants with respect to the
Chevalley basis are being reduced modulo 3: the simple pieces of
the superalgebras thus obtained appear on Elduque's super square
as $\fg(2, 3)$ and $\fg(2, 6)$ in notations described below. Just
compare the Cartan matrices 3) of $\fa\fg(2)$ and 6) of
$\fa\fb(3)$ in \cite{GL1} with the Cartan matrices 3) of $\fg(2,
3)$ and 6) of $\fg(2, 6)$, respectively: They are identical.
Observe the crucial difference between the Lie (super)algebra
$\fg(A)_\Kee$ (which might be an algebra without Cartan matrix)
and the Lie (super)algebra $\fg(A_p)$ constructed from the integer
Cartan matrix $A$ reduced modulo $p$ if $p<5$. In particular, for
any field $\Kee$ of characteristic $3$, the superdimensions (and
structures) of $\fg(2, 3)$ and $\fa\fg(2)_\Kee$ are distinct as
well as those of $\fg(2, 6)$ and $\fa\fb(3)_\Kee$.

Elduque and Cunha gave only one Cartan matrix for each Lie
superalgebra that possess Cartan matrix, and here we list all
possible inequivalent Cartan matrices (in other words,
inequivalent systems of simple roots). Among all inequivalent
systems of simple roots, some have properties particularly useful
in certain applications (e.g., the two extremes: with the least
number of odd roots and with the maximal number of odd roots).

We also need the list of all systems of simple roots as a step
towards classification of prolongations of the non-positive
components of the Lie superalgebras we consider here in their
various gradings (for details of this approach, see \cite{GL4}).

To interpret these new Elduque-and-Cunha Lie superalgebras as Lie
superalgebras preserving a tensor is still a task to be performed.
This interpretation is not, however, always possible; in the
sequel to this paper we intend to interpret these Lie
superalgebras as preserving not just a tensor but certain
distributions, as in \cite{GL4}, and get new examples of simple
Lie superalgebras, as in \cite{BjL}. Here, with the aid of the
{\bf SuperLie} package, we describe presentations of those
Elduque-and-Cunha Lie superalgebras that possess Cartan matrix.

In terms of the Chevalley generators of the Lie superalgebras with
Cartan matrix, there are two types of defining relations:
Serre-type ones and non-Serre type ones (over $\Cee$, all the
relations are listed in \cite{GL1}). Sometimes some of the
Serre-type relations are redundant but this does not matter in
practical calculations. At the moment, the problem \so{how to
encode the non-Serre type relations in terms of Cartan matrix} is
open. Some relations are so complicated that we conjecture that
there is no {\it general} encoding procedure. This is why our list
of relations is of practical interest.

In the \textrm{arXiv} version of this paper, we also give the
coefficients of linear dependence over $\Zee$ of the corresponding
maximal roots with respect to the simple ones and the inverses of
Cartan matrices whenever the latter are invertible.

\section{Background}

\ssec{2.1}{Notations} The ground field $\Kee$ of characteristic
$p$ is assumed to be algebraically closed. We follow Bourbaki's
convention: if $G$ is a Lie group, then its Lie algebra is
designated $\fg$, although the modern tradition does not favor
Gothic font  in characteristic $p>0$. Let $\fg':=[\fg, \fg]$.

\ssec{2.2}{What $\fg(A)$ is} First, recall, how to construct a Lie
superalgebra from a Cartan matrix (\cite{Se,vdL,FLS,Eg,FSS}). Let
$A=(A_{ij})$ be an arbitrary $n \times n$ matrix of rank $l$ with
entries in $\Kee$. Fix a vector space $\fh$ of dimension $2n-l$
and its dual $\fh^{*}$, select $n$ linearly independent vectors
$h_{i} \in \fh$ and let $\alpha_{j} \in \fh^{*}$ be such that
$\alpha_{i}(h_{j}) = A_{ij}$.

Let $I= \{i_{1},\dots, i_{n}\} \subset (\Zee/2\Zee)^{n}$; consider
the free Lie superalgebra $ \widetilde{\fg}(A, I)$ with generators
$e_{1}^\pm,\dots , e_{n}^\pm$, where $p(e_{j}^\pm)=i_{j}$, and
defining relations (for convenience, we set $h_i:=[e_{i}^+,
e_{i}^-]$):
\begin{equation}
\label{relcm} [e_{i}^+, e_{j}^-] = \delta_{ij} h_{j};\quad [h_{i},
e_{j}^\pm]=\pm A_{ij}e_{j};\quad \;[h_i, h_j]=0.
\end{equation}

The following statements over $\Cee$ is well known for the Lie
algebras; for Lie superalgebras over $\Cee$ it is due to Serganova
and van de Leur \cite{Se, vdL}; for for Lie superalgebras over
$p\neq 2$, see \cite{CE2}. (The proof does not seem to depend on
the super structure or $p$ if $p\neq 2$.)

\begin{Statement} {\em a)} Let $\tilde\fn_{+}$
and $\tilde\fn_{-}$ be the superalgebras in $\tilde\fg(A, I)$
generated by $ e_{1}^\pm$,\dots , $e_{n}^\pm$; then $\tilde\fg(A,
I) \cong \tilde\fn_{+} \oplus \fh\oplus \tilde\fn_{-}$, as vector
superspaces.

{\em b)} Among the ideals of $\tilde\fg (A, I)$ whose intersection
with $\fh$ is zero, there exists a maximal ideal $\fr$ such that
$\fr$ is the direct sum of ideals $\fr\bigcap\tilde\fn_{+}$ and
$\fr\bigcap\tilde\fn_{-}$.
\end{Statement}

Set $\fg(A, I) =\tilde\fg (A, I)/\fr$. Both $\fg(A, I)$ and
$\fg'(A, I)$ may contain a center. As proved in \cite{Se,vdL}, the
centers $\fc$ of $ \fg(A, I)$ and $\fc'$ of $ \fg'(A, I)$ consist
of all $h \in \fh$ such that $\alpha_{i}(h) =0$ for all $i=
1,\dots , n$; this is also true for $p>0$.

Clearly,
\begin{equation}
\label{rescale} \text{the rescaling
$e_i^\pm\mapsto\sqrt{\lambda_i}e_i^\pm$, sends $A$ to
$\diag(\lambda_1, \dots , \lambda_n)\cdot A$.}
\end{equation}
Two pairs $(A, I)$ and $(A', I')$ are said to be {\it equivalent}
if $(A', I')$ is obtained from $(A, I)$ by a composition of a
permutation of indices and a rescaling $A' = \diag (\lambda_{1},
\dots, \lambda_{n})\cdot A$, where $\lambda_{1}\dots
\lambda_{n}\neq 0$. Clearly, equivalent pairs determine isomorphic
Lie superalgebras.

The matrix $A$ (more precisely, the pair $(A, I)$) is said to be a
{\it Cartan matrix} of the Lie superalgebra $\fg (A, I) $ and also
of $\tilde\fg(A, I)$, $\fg'(A, I)$, as well as of $\fg(A, I)/\fc$
and $\fg' (A, I)/ \fc'$. We tirelessly repeat:

\so{The} \lq\lq \so{relatives}\rq\rq\ (\so{such as nontrivial central
extensions and algebras of derivations}) \so{of simple Lie}
(\so{super})\so{algebras are no less important than the simple algebras
themselves}.

In particular, in Elduque's superization of the Freudenthal Magic
Square for $p=3$, there naturally appear $\fg(2, 8)=\fe(6)$,
$\fg(2, 3)$ and $\fg(2, 6)$ that have 1-dimensional centers.

\ssec{2.3}{Systems of simple roots} Since $\fg(A, I)$, where $A$
is of size $n\times n$, is naturally $\Zee^n$-graded, it follows
that we can consider the root system $R$ of $\fg(A, I)$ as a
subset of $\Zee^n$.  For any subset $B=\{\sigma_{1}, \dots,
\sigma_{n}\} \subset R$, we set (we denote by $\Zee_{+}$ the set
of non-negative integers):
$$
R_{B}^{\pm} =\{ \alpha \in R \mid \alpha = \pm \sum n_{i}
\sigma_{i},\;\;n_{i} \in \Zee_{+} \}.
$$

The set $B$ is called a {\it system of simple roots} \index{root,
simple, system of} of $R$ (or $\fg$) if $ \sigma_{1}, \dots ,
\sigma_{n} $ are linearly independent and there exist $\tilde
e_{1}^\pm\,{\in}\;\,\fg _{\pm\sigma_{1}},\dots , \tilde
e_{n}^\pm\,{\in}\;\,\fg_{\pm\sigma_{n}}$ such that:
$$
\fg = \fg _{B}^{-} \oplus \fh \oplus \fg _{B}^{+},
$$
where $\fg _{B}^{\pm}$ is the super\-algebra gen\-er\-ated by
$\tilde e_{1}^\pm,\dots , \tilde e_{n}^\pm$.

Let $B$ be a system of simple roots and $\tilde e_{1}^\pm, \dots ,
\tilde e_{n}^\pm$ the corresponding elements of $\fg(A, I)$.  Set
$\tilde h_{i}=[\tilde e_{i}^{+}, \tilde e_{i}^-]$ , $A_{B}
=(A_{ij})$, where $A_{ij} =\sigma_{i}(\tilde{h_{j}})$ and
$I_{B}=\{p(\sigma_{1}), \cdots, p(\sigma_{n})\}$. The Cartan
matrix $(A_{B}, I_{B})$ thus constructed does not have to coincide
with the initial Cartan matrix $(A, I)$.


Two systems of simple roots $B_{1}$ and $B_{2}$ are said to be
{\it equivalent} if the pairs $(A_{B_{1}}, I_{B_{1}})$ and
$(A_{B_{2}}, I_{B_{2}})$ are equivalent.

Hereafter $\fg =\fg (A, I)$. To be able to distinguish the case of
the even root from the odd one when $A_{ii}=0$, we write
$A_{ii}=\emptyset$ if $p(\sigma_{i})=\ev$. The following statement
is subject to a direct verification:

\begin{Proposition} Let $B$ be a system of simple roots of $\fg$, $\tilde e_{i}^\pm$
for $i=1$,\dots, $n$ the corresponding set of generators and
$A_{B}=(A_{ij})$ the Cartan matrix. Fix an $i$. Then:

{\em a)} If $p(\sigma_{i})= \ev $ and $A_{ii} \neq 0$, then the
Lie subalgebra generated by the $e_{i}^\pm$ is isomorphic to
$\fsl(2)$.

{\em b)} If $p(\sigma_{i})= \ev $ and $A_{ii} = \emptyset$, then
the Lie subalgebra generated by the $e_{i}^\pm$ is isomorphic an
analog of the Heisenberg Lie algebra; we denote this Lie algebra
$\fhei(2; p; \underline{N})$, where $N\in \Nee$. Its natural
representation is realized in the Fock space of functions $\cO(1;
\underline{N})$; it is indecomposable\footnote{At least, under the
conventional definition of representations of Lie algebras for
$p>0$; for conjectural other versions, see \cite{LL}.} for
$\underline{N}>1$ and irreducible for $\underline{N}=1$.

{\em c)} If $p(\sigma_{i}) =\od $ and $A_{ii}=0$, then
$2\sigma_{i} \not \in R$ and the subsuperalgebra generated by the
$e_{i}^\pm$ is isomorphic to $\fsl(1|1)$.

{\em d)} If $p(\sigma_{i}) =\od $ and $A_{ii} \neq 0$, then
$3\sigma_{i} \not \in R$ and the subsuperalgebra generated by the
$e_{i}^\pm$ is isomorphic to $\fosp(1|2)$.
\end{Proposition}
\ssec{2.4}{Chevalley generators and odd reflections} How many
Cartan matrices correspond to the same Lie superalgebra $\fg$? To
answer this question, let us multiply $A_{B}$ by a diagonal matrix
as in $(\ref{rescale})$, so that in the cases a), d) of
Proposition 2.3 the diagonal elements of $A_{B}$ become 2 or 1,
respectively. The row with a 0 on the main diagonal can be
multiplied by any nonzero factor; we usually multiply it so that
$A_{B}$ be symmetric  if possible, and additionally, if $p=0$, so
as to make the off-diagonal elements negative (this is needed in
relations~$(\ref{srpm})$). Such a matrix is said to be {\it
normalized}.

For each simple finite dimensional Lie algebras (and,
conjecturally, for Lie superalgebras) for $p>3$ the normalized
Cartan matrix might be not symmetric, but it is symmetrizable. We
only consider normalized matrices, unless we need them in a
symmetric form. A symmetrized matrix is not a Cartan matrix but is
also useful: it gives the values of the inner products of simple
roots and is needed to pass from one system of simple roots to
another.

A usual way to represent simple Lie algebras over $\Cee$ with
integer Cartan matrices is via graphs called, in the finite
dimensional case, {\it Dynkin diagrams}.  The Cartan matrices of
Lie {\it super}algebras $\fg$ (even over $\Cee$) can be
non-symmetrizable or (for any $p$ in the super case and for $p>0$
in the  non-super case) have entries belonging to the ground field
$\Kee$. Still, it is always possible to assign an analog of the
Dynkin diagram to both Lie superalgebras of polynomial growth, and
to finite dimensional Lie algebras (if these (super)algebras
possess Cartan matrices). Perhaps, the edges and nodes of the
graph should be rigged with an extra information. Namely, the {\it
Dynkin--Kac diagram}\index{Dynkin--Kac diagram} of the matrix $(A,
I)$ is a set of $n$ nodes connected by multiple edges, perhaps
endowed with an arrow, according to the usual rules (\cite{K}) or
their modification, most naturally formulated by Serganova: cf.
\cite{Se, FLS} with \cite{FSS}.

The nodes are of four types: {\sl To every simple root there
corresponds}
$$
\begin{cases}
\text{a node}\; \motimes \;& \text{ if }\; p(\alpha_{i}) =\od \;
\text{ and }\;
A_{ii}=0,\\
\text{a  node}\; \mbullet \;&\text{ if }\;
p(\alpha_{i}) =\od \; \text{ and }\; A_{ii}=1;\\
\text{a  node}\; \mcirc\; &\text{ if }\; p(\alpha_{i})= \ev ,\;
\text{ and }\;
A_{ii}=2\\
\text{a  node}\; *\; &\text{ if }\; p(\alpha_{i})= \ev ,\; \text{
and }\;
A_{ii}=\emptyset.\\
\end{cases}
$$
{\it A posteriori} (from the classification of simple Lie
superalgebras with Cartan matrix and of polynomial growth for
$p=0$) we find out that the roots~$*$ can only occur if $\fg(A,
I)$ grows faster than polynomially. Thanks to classification
again, if $\dim \fg<\infty$, the roots $*$ can not occur if $p>3$;
whereas for $p=3$, the Brown Lie algebras are examples of $\fg(A)$
with a simple root of type roots $*$, see \cite{GL4}.

To more graphically express normalized Cartan matrices, we apply
Serganova's rules even for $p>0$. Although the analog of the
Dynkin graph is now uniquely recovered from the Cartan matrix only
if there are no two odd \lq\lq grey" nodes~$\motimes$ connected,
this analog helps to grasp the geometry of the system of simple
roots.

We often denote the set of generators corresponding to a
normalized matrix by $X_{1}^{\pm},\dots , X_{n}^{\pm}$ instead of
$e_{1}^{\pm},\dots , e_{n}^{\pm}$; and call them the {\it
Chevalley generators}.

Let $\alpha$ be a root of $\fg(A, I)$ with the entries of $A$ in
$\Zee$ or $\Zee/p$. Define the reflections $r_{\alpha}$ acting on
the root lattice over $\Zee$ by the formulas
\begin{equation}
\label{oddrefl} r_{\alpha_{i}}(\alpha_{j})= \left\{
\begin{array}{lllll}
-\alpha_{j} &\text{for $i=j$}\\
\alpha_{j}-A_{ij}\alpha_{i} &\text{for $i \neq j$ and $A_{ii}
=2$}\\
\alpha_{j}-2A_{ij}\alpha_{i} & \text{for $i \neq j$ and $A_{ii}=1$} \\
\alpha_{j}+\alpha_{i} & \text{for $i \neq j$ and $A_{ii}=0
\;,\;A_{ji} \neq
0$}\\
\alpha_{j} &\text{for $i \neq j$ and $A_{ii}=A_{ji}=0$}
\end{array} \right.
\end{equation}

The reflection in the $i$th root sends one set of Chevalley
generators into the new one:
\begin{equation}
\label{oddrefx} \tilde X_{i}^{\pm}=X_{i}^{\mp};\;\;
\tilde X_{j}^{\pm}=\begin{cases}[X_{i}^{\pm}, X_{j}^{\pm}]&\text{if $A_{ij}\neq 0$}\\
X_{j}^{\pm}&\text{otherwise}.\end{cases}
\end{equation}

The reflections in the odd roots are referred to as {\it odd
reflections}. For $p=0$, {\sl there is always a chain of odd
reflections connecting every two systems of simple roots $B_1$ and
$B_2$}. Serganova's proof of this fact \cite{Se} is applicable
without modifications for $p\neq 2$. Theorem 1.5.5 in \cite{S}
ensures that the maximal tori of $\fg_\ev$ are conjugate, so by
means of odd reflections we do get all inequivalent Cartan
matrices.

\ssec{2.5}{Defining relations for Lie superalgebras $\fg=\fg(A,
I)$} The simple Lie superalgebras of the form $\fg=\fg(A, I)$ have
several quite distinct sets of relations (cf. \cite{Sa} and refs
therein) but usually they are given by their {\it Chevalley
generators} $X^{\pm}_i$ of degree $\pm 1$ to which the elements
$H_i = [X_{i}^+, X_{i}^-]$ are added for convenience. These
generators satisfy the following relations (hereafter in similar
occasions either all superscripts $\pm$ are $+$ or all are $-$)
\begin{equation}
\label{g(a)} {}[X_{i}^+, X_{j}^-]  = \delta_{ij}H_i, \quad [H_{i},
H_{j}] = 0, \quad [H_i, X_j^{\pm}] = \pm A_{ij}X_j^{\pm},
\end{equation}
and additional relations $R_i=0$ whose left sides are implicitly
described, for a general Cartan matrix with entries in $\Kee$, as
(\cite{K})
\begin{equation}
\label{myst}
\begin{split}
 &\text{\lq\lq the $R_i$ that generate the maximal ideal $I$ such that}\\
 &I\cap\Span(H_i \mid 1\leq i\leq n)=0. \text{\rq\rq}
\end{split}
\end{equation}
For $p=0$ and normalized Cartan matrices of simple finite
dimensional Lie algebras, there is only one ({\it Chevalley})
basis in which all structure constants are integer, cf.
\cite{Er}.(Having normalized the Cartan matrix of $\fo(2n+1)$ so
that $A_{nn}=1$ we get another basis with integer structure
constants.)

We conjecture that if $p>2$, then for Lie algebra and Lie
superalgebras with Cartan matrix with entries in $\Zee/p$, there
also exists only one analog of the Chevalley basis. At least, one
such basis definitely exists: we use it.

\subsubsection{Serre relations} Let $\fg=\fg(A, I)$; let
$X_{i}^{\pm}$, where $1\leq i \leq n$, be root vectors
corresponding to simple roots $\pm\sigma_i$. Clearly, the
$X_{i}^{\pm}$ generate $\fn^{\pm}$. We find the defining relations
by induction on $n$ with the help of the Hochschield--Serre
spectral sequence (for its description for Lie superalgebras,
which has certain subtleties, see \cite{Po}). For the basis of
induction consider the following cases:
\begin{equation}
\label{3.4.1}
\begin{array}{lll}
   \mcirc \;\; or\;\; \mbullet &\text{no relations, i.e., $\fn^{\pm}$
are free Lie superalgebras}&\text{if $p\neq 3$;}\\
\mbullet &\left(\ad_{X^{\pm}}\right)^2(X^{\pm})=0&\text{if $p=3$;}\\
{\motimes}& [X^{\pm}, X^{\pm}]=0.&
\end{array}
\end{equation}
Set $\deg X_{i}^{\pm} = 0$ for $ 1\leq i\leq n-1$ and $\deg
X_{n}^{\pm} = \pm 1$.  Let $ \fn^{\pm} = \oplus \fn_{ i}^{\pm}$
and $\fg= \oplus \fg_{i} $ be the corresponding $\Zee $-gradings.
Set $\fn_{\pm} =\fn^{\pm}/\fn_{ 0}^{\pm}$.  From the
Hochschield--Serre spectral sequence for the pair $\fn_{ 0}^{\pm}
\subset \fn ^{\pm}$ we get:
\begin{equation}
\label{3.4.2}
H_{2}(\fn_{0}^{\pm})\oplus H_{1}(\fn_{ 0}^{\pm}; H_{
1}(\fn_{\pm}))\oplus H_{0}(\fn_ { 0}^{\pm}; H_{ 2}(\fn_{\pm})).
\end{equation}
It is clear that
\begin{equation}
\label{3.4.3}
H_{1}(\fn_{\pm})= \fn_{ 1}^{\pm} , \;\;\; H_{
2}(\fn_{\pm}) = E^{2}(\fn_ { 1}^{\pm})/\fn_{ 2}^{\pm},
\end{equation}
where $E^2$ is the functor of the second exterior power. So, the
second summand in~(\ref{3.4.2}) provides us with relations of the
form:
\begin{equation}
\label{3.4.4} \begin{array}{ll} (\ad_{X_{n}^{\pm}})^{k_{ni}}
(X_{i}^{\pm})=0&\text{if the $n$-th root is not}\;\;
{\motimes}\\
{} [X_{n}, X_{n}]=0&\text{if the $n$-th root is} \;\;\motimes.
\end{array}
\end{equation}
while the third summand in (\ref{3.4.2}) is spanned by the
$\fn_{0}^{\pm}$-lowest vectors in
\begin{equation}
\label{3.4.6} E^{2}(\fn_{1}^{\pm})/(\fn_{ 2}^{\pm} + \fn^{\pm}
E^{2}(\fn_{1}^{\pm})).
\end{equation}

Let the matrix $B=(b_{ij})$ be obtained from the Cartan matrix
$A=(A_{ij})$ by replacing all nonzero elements in the row with
$A_{ii}=0$ by $-1$ and multiplying the row with $A_{ii}=1$ by $2$.
The following proposition, whose proof is straightforward,
illustrates the usefulness of our normalization of Cartan matrices
as compared with other options:

\begin{Proposition} The numbers $k_{in}$ and $k_{ni}$ in $(\ref{3.4.4})$ are
expressed in terms of $(b_{ij})$ as follows:
\begin{equation}
\label{srpm}
\begin{matrix} (\ad_{X_{i}^{\pm}})^{1-b_{ij}}(X_{j}
^{\pm})=0 &
\text{ for $ i \neq j$} \\ & \\
\text{$[X_{i}^{\pm}, X_{i}^{\pm}]=0$} & \text{if $A_{ii}
=0$}.\end{matrix}
\end{equation}
\end{Proposition}
The relations $(\ref{g(a)})$ and $(\ref{srpm})$ will be called {\it
Serre relations} for Lie superalgebra $\fg(A, I)$. If $p=3$, then
the relation
\begin{equation}
\label{x^3}
[X_{i}^{\pm}, [X_{i}^{\pm}, X_{i}^{\pm}]]=0
\;\;\text{for $X_{i}^{\pm}$ odd and $A_{ii} =1$}
\end{equation}
is not a consequence of the Jacobi identity; for simplicity, we
will incorporate it to the set of Serre relations.

\subsubsection{Non-Serre relations} correspond to the
third summand in~(\ref{3.4.2}).  Let us consider the simplest
case: $\fsl (m|n)$ in the realization with the system of simple
roots
\begin{equation}
\label{circ}
 \xymatrix@C=1em{
 \ffcirc\ar@{-}[r]&\dots\ar@{-}[r]&\ffcirc\ar@{-}[r]&
 \otimes\ar@{-}[r]&\ffcirc\ar@{-}[r]&\dots\ar@{-}[r]&
 \ffcirc
 }
\end{equation}

Then $H_2(\fn_{\pm})$ from the third summand in (\ref{3.4.2}) is
just $E^2(\fn_{\pm})$.  For simplicity, confine ourselves to the
positive roots.  Let $X_{1}$, \dots , $X_{m-1}$ and
$Y_{1}$,~\dots,~$Y_{n-1}$ be the root vectors corresponding to
even roots separated by the root vector $Z$ corresponding to the
root $\motimes$.

If $n=1$ or $m=1$, then $E^2(\fn)$ is an irreducible $\fn_{\bar
0}$-module and there are no non-Serre relations.  If $n\neq 1$ and
$m\neq 1$, then $E^2(\fn)$ splits into 2 irreducible $\fn_{\bar
0}$-modules.  The lowest component of one of them corresponds to
the relation $[Z, Z]=0$, the other one corresponds to the
non-Serre-type relation
\begin{equation}
\label{*} [[X_{m-1}, Z], [Y_{1}, Z]] =0. \end{equation}

If, instead of $\fsl (m|n)$, we would have considered the Lie
algebra $\fsl(m+n)$, the same argument would have led us to the
two relations,  both of Serre type:
$$
\ad_Z^2(X_{m-1})=0, \qquad \ad_Z^2(Y_{1})=0.
$$
Although we have found all the {\it basic} relations, the ones
that generate the ideal of relations, in what follows, we only
list non-Serre relations assuming that \so{all} Serre relation are
satisfied despite the fact that some of the Serre relations turn
out to be consequences of these basic relations. Such redundances
are rare, to single them out is a boring task of doubtful value.

\section{Elduque and Cunha superalgebras: Systems of simple roots}

For details of description of Elduque and Cunha superalgebras in
terms of symmetric composition algebras, see \cite{El1,CE,CE2}.
Here we consider the simple Elduque and Cunha superalgebras with
Cartan matrix for $p=3$. In what follows, we list them using
somewhat shorter notations as compared with the original ones:
here $\fg(A,B)$ denotes the superalgebra occupying $(A, B)$th slot
in the Elduque--Freudenthal magic super square; the first Cartan
matrix is usually the one given in \cite{CE}, the other matrices
are obtained from the first matrix by means of odd reflections.
Accordingly, $\fg(A, B)^{c)}$ is the shorthand for the realization
of $\fg(A, B)$ by means of Cartan matrix $c)$.

The table at the beginning of each case shows the result of odd
reflections (the number of the row is the number of the matrix in
the list below, the number of the column is the the number of the
root in which reflection is made; the cells contain the results of
reflections (the number of the matrix obtained) or a \lq\lq --" if
the reflection is not appropriate because $A_{ii}\neq 0$. Some of
the matrices thus obtained are equivalent; we did not eliminate
a few redundancies.

If the diagram of $\fg$ is symmetric, it gives rise to an outer
automorphism whose fixed points constitute a Lie superalgebra. The
examples below where this occurred did not lead to any new simple
Lie superalgebra.

The numbers of matrices with the maximal number of even roots are
boxed, those with the maximal number of add roots are underlined.
The nodes are numbered by small boxed numbers; the dashed lines
with arrows depict odd reflections.

\begin{figure}[ht]\centering
\parbox{.28\linewidth}{$
\fg{}(2,3)\;\; \begin{lmatrix}
    -&-&2 \\
    3&4&1 \\
    2&5&- \\
    5&2&- \\
    4&3&-
  \end{lmatrix}
$}\hfill
\parbox{.7\linewidth}{\includegraphics{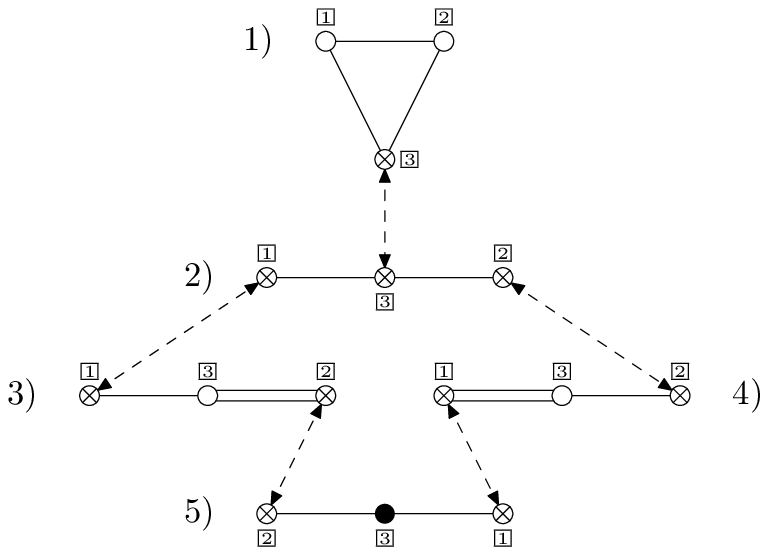}}

\end{figure}

{\footnotesize
$$

$$
}


\begin{figure}[ht]\centering
\parbox{.5\linewidth}{\mbox{}\hfill$
\fg{}(1,6)\;\;
\begin{lmatrix}
    -&-&2 \\
    -&-&1
  \end{lmatrix}
$\quad\mbox{}}\hfill
\parbox{.5\linewidth}{\includegraphics{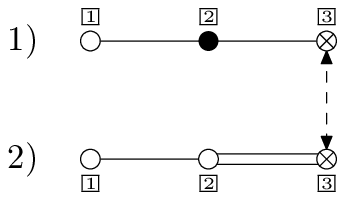}}
\end{figure}

\vskip-4mm

{\footnotesize
$$ 1) %
$$

\begin{figure}[ht]\centering
\includegraphics{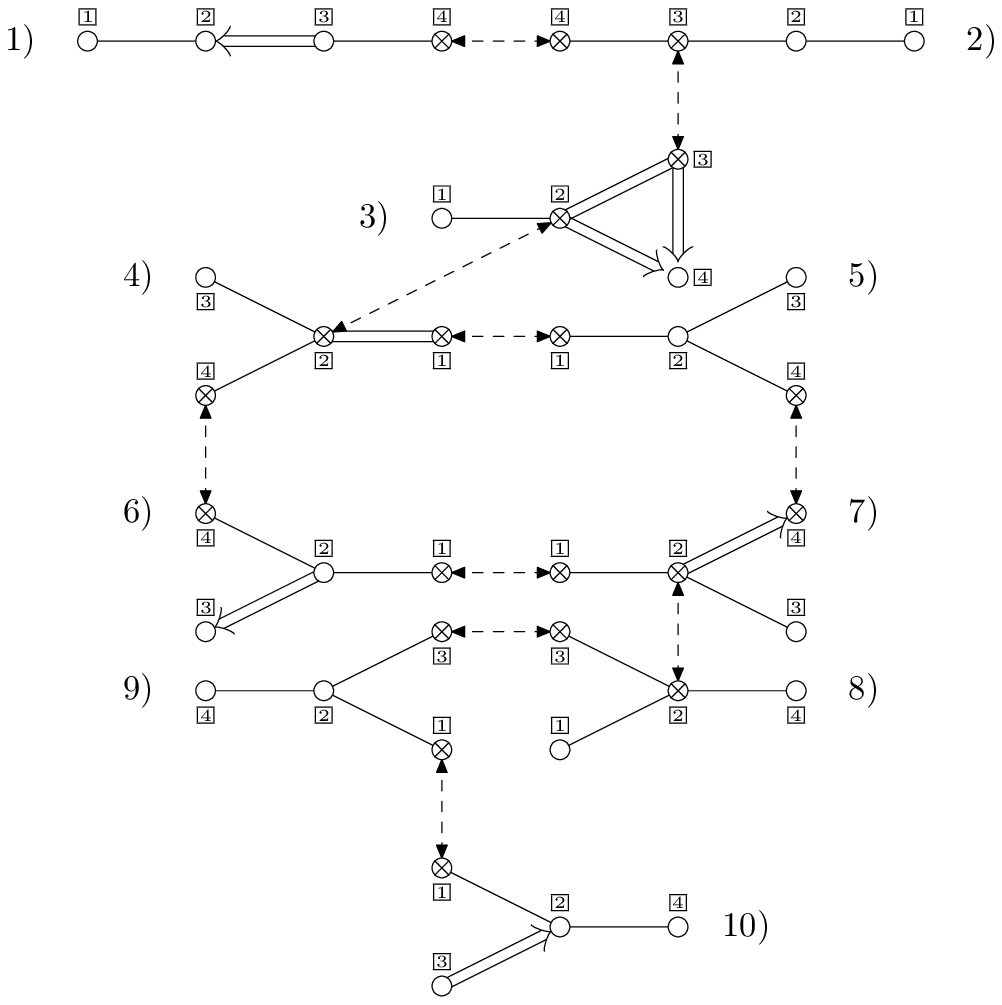}
\end{figure}

\tiny

$$
%
  }
$$

\vskip-3mm

\begin{figure}[ht]\centering
\includegraphics{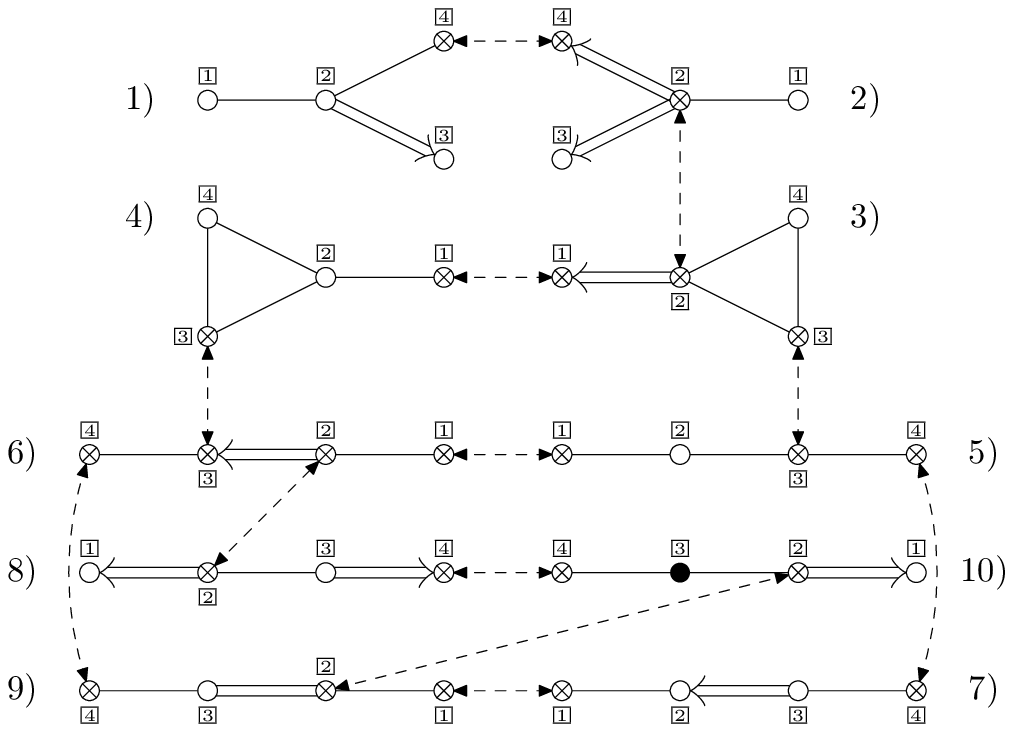}
\end{figure}

\vskip-3mm

\tiny
$$
%
  }
  $$

\begin{figure}[ht]\centering
\includegraphics{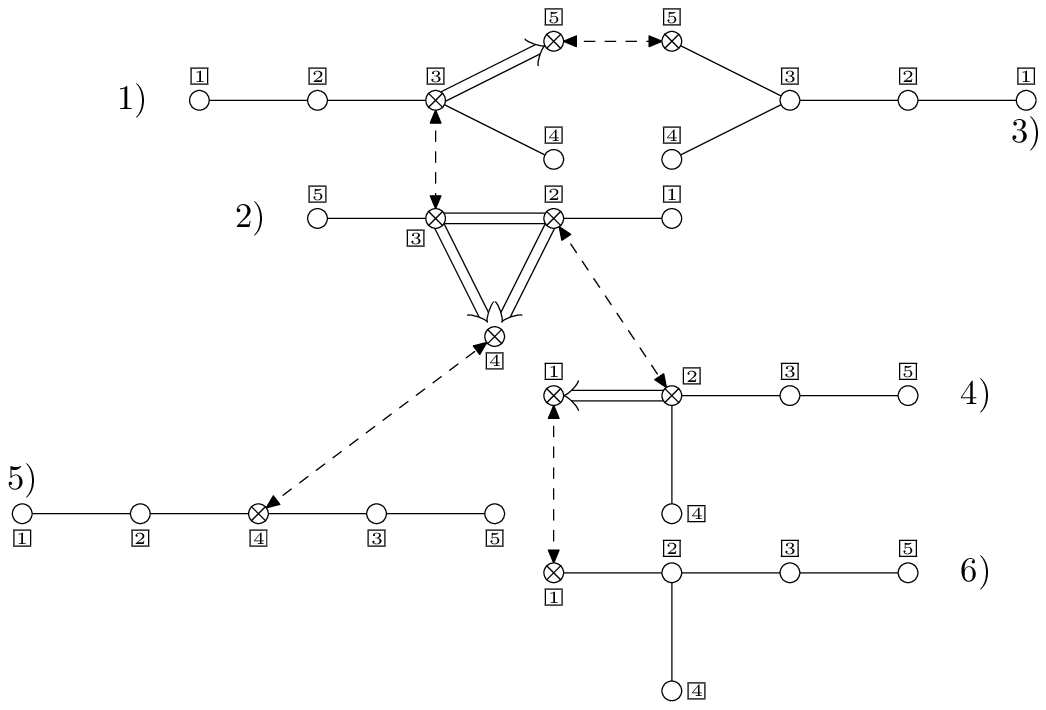}
\end{figure}

\tiny

$$%
$$

\begin{figure}[ht]\centering
\includegraphics{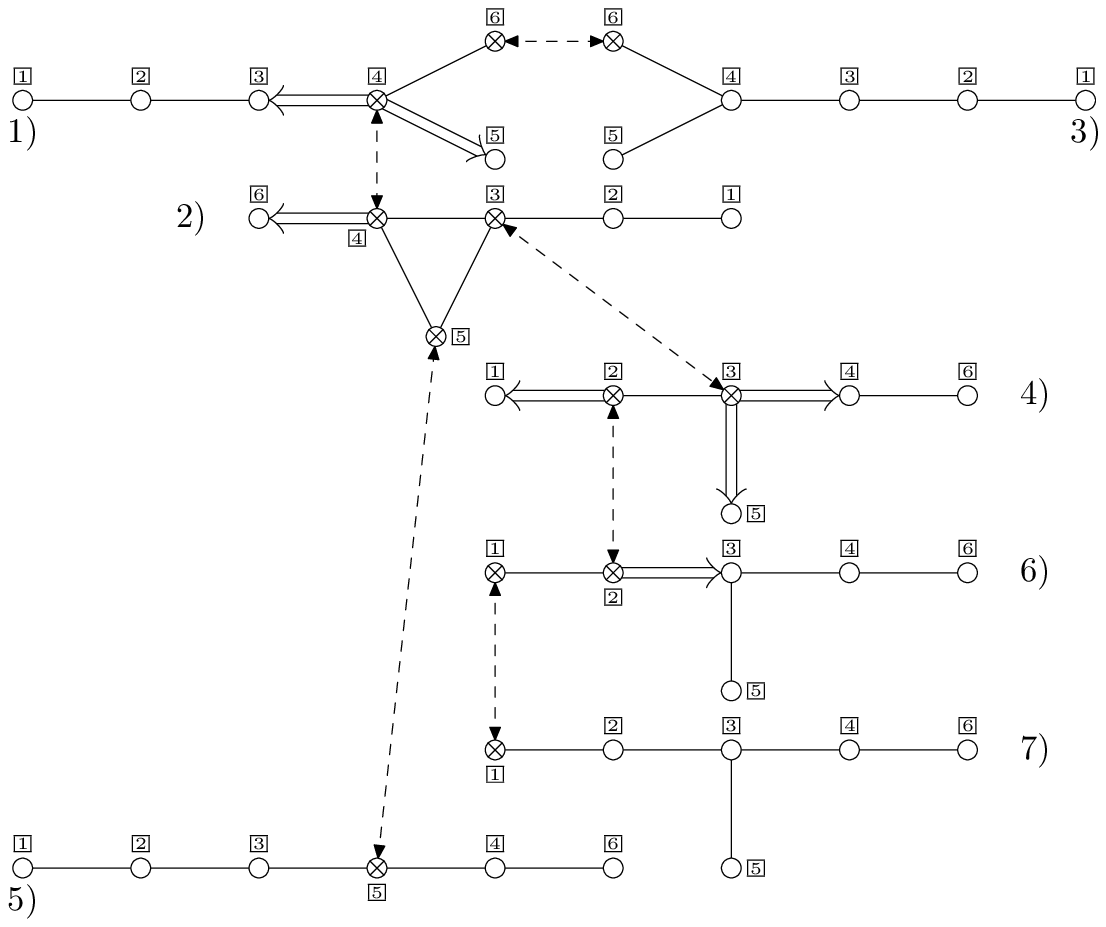}
\end{figure}

\vspace*{-5mm}

\tiny
$$
%
$$

\begin{figure}[ht]\centering
\includegraphics{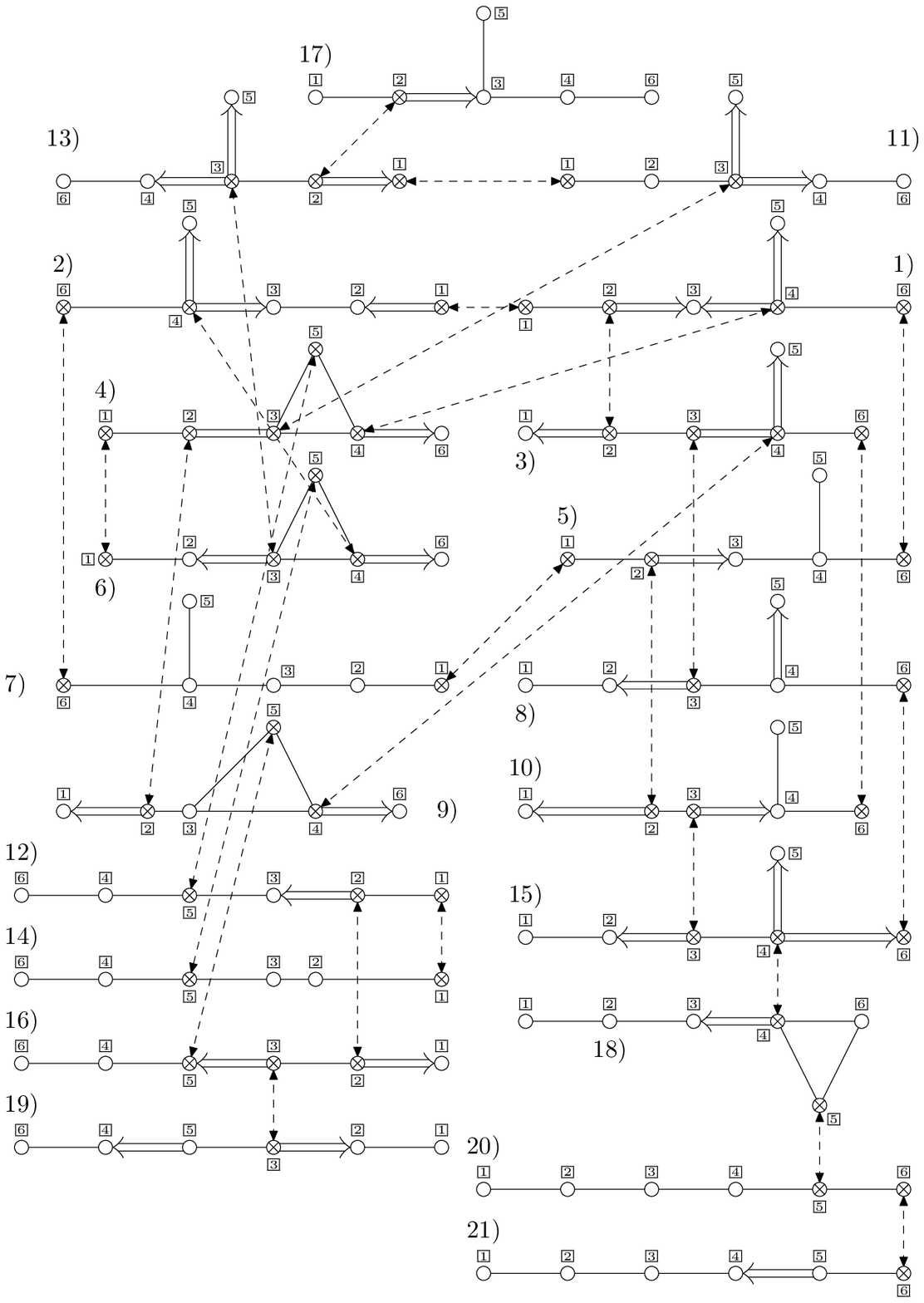}
\end{figure}

\newpage
\vspace*{-\baselineskip}

\tiny

$$
\fg{}(8,6) \quad %
$$

\begin{landscape}
\begin{figure}[ht]\centering
\includegraphics{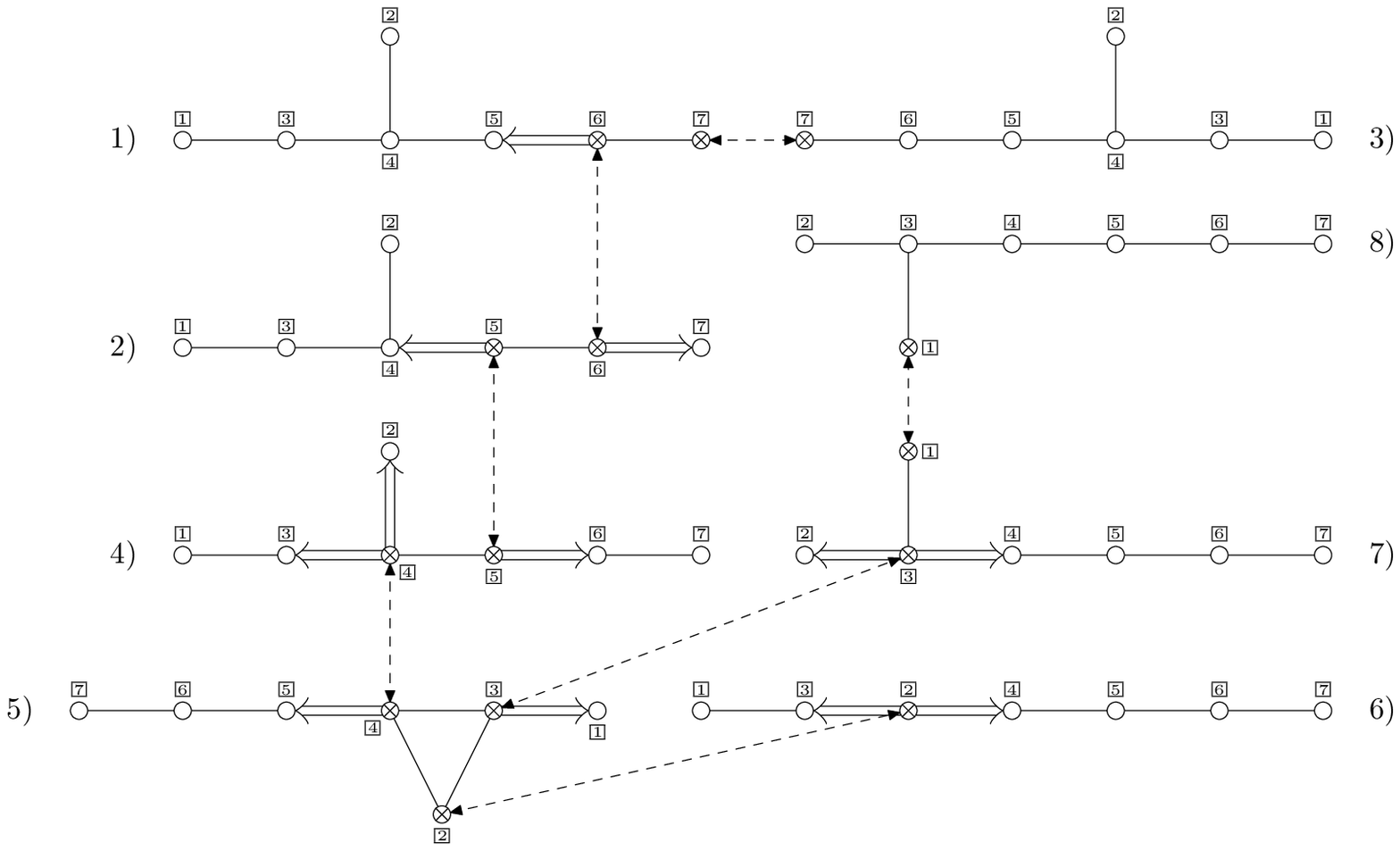}
\end{figure}
\end{landscape}


\medskip

\normalsize

\section{Elduque and Cunha superalgebras: Defining relations}
To save space, in what
follows we omit indicating the Serre relations; their fulfilment is assumed. 

\footnotesize{
\begin{arab}

\item[]

\hskip2.5ex\underline{$\fg(2, 3), \;\sdim =(11|14)$}

\item[]
$
\arraycolsep=0pt
\renewcommand{\arraystretch}{1.2}
\begin{array}{ll}
1)\ &{{[x_{2},\,[x_{1},\,x_{3}]]}-{2\,
[x_{3},\,[x_{1},\,x_{2}]]}}=0;\\ &
\left(\ad_{[x_{1},\,x_{2}]}\right)^2([x_{1},\,x_{3}])=0;\quad
\left(\ad_{[x_{1},\,x_{2}]}\right)^2([x_{2},\,x_{3}])=0.\\
      \end{array}
$

$
\arraycolsep=0pt
\renewcommand{\arraystretch}{1.2}
\begin{array}{ll}
2)\ &   \left (\ad_{[x_1,x_3] } \right )^2([x_2,x_3])=0;\quad
 \left(\ad_{[x_2,x_3]} \right )^2([x_1,x_3])=0;\\
& {{[[x_{2},\,[x_{1},\,x_{3}]],\,
      [[x_{1},\,x_{3}],\,[x_{2},
      \,[x_{1},\,x_{3}]]]]}= 0}; \\
&  {{[[x_{2},\,[x_{1},\,x_{3}]],\,
      [[x_{2},\,x_{3}],\,[x_{2},\,
      [x_{1},\,x_{3}]]]]}=0.}\\
      \end{array}
$

$
\arraycolsep=0pt
\renewcommand{\arraystretch}{1.2}
\begin{array}{ll}
3)\ &
{{[[x_{1},\,x_{3}],\,[x_{3},\,[x_{2},\,x_{3}]]]}=0};\quad
  {\left ( \mathrm{ad}_{[x_2,x_3]}\right )^2(x_3) =0}; \\
&  {{[[x_{2},\,[x_{1},\,x_{3}]]\!,\,
    [[x_{1},\,x_{3}]\!,\,[x_{2},\,x_{3}]]]}\,{=}\,0};\\
&  {{[[x_{2},\,[x_{1},\,x_{3}]]\!,\,
    [[x_{2},\,x_{3}]\!,\,[x_{2},\,x_{3}]]]}\,{=}\,0.}
\end{array}
$

$
\arraycolsep=0pt
\renewcommand{\arraystretch}{1.2}
\begin{array}{ll}
4) \ &  \left ( \mathrm{ad}_{[x_1,x_3]}\right )^2(x_3)=0; \quad
  {{[[x_{2},\,x_{3}],\,[x_{3},\,[x_{1},\,x_{3}]]]}=0};\\
&  {{[[x_{2},\,[x_{1},\,x_{3}]]\!,\,
    [[x_{1},\,x_{3}]\!,\,[x_{1},\,x_{3}]]]}\,{=}\,0};\\
&  {{[[x_{2},\,[x_{1},\,x_{3}]]\!,\,
    [[x_{1},\,x_{3}]\!,\,[x_{2},\,x_{3}]]]}\,{=}\,0.}
\end{array}
 $

 $
\arraycolsep=0pt
\renewcommand{\arraystretch}{1.2}
\begin{array}{ll}
5)\ & {{[x_{3},\,[x_{2},\,[x_{1},\,x_{3}]]]}-
    {[[x_{1},\,x_{3}],\,[x_{2},\,x_{3}]]}=0};\\
&
{{[[x_{1},\,x_{3}],\,[x_{3},\,x_{3}]]}=0};\quad
  {{[[x_{2},\,x_{3}],\,[x_{3},\,x_{3}]]}=0};\\
& {{[[x_{2},\,[x_{1},\,x_{3}]],\,[x_{3},\,
  [x_{1},\,x_{3}]]]}=0};\quad
  {{[[x_{2},\,[x_{1},\,x_{3}]],\,
  [x_{3},\,[x_{2},\,x_{3}]]]}=0.}
    \end{array}
$

\item[]

\hskip2.5ex \underline{$\fg(1,6), \;\sdim =(21|14)$}

\item[]

$
\arraycolsep=0pt
\renewcommand{\arraystretch}{1.2}
\begin{array}{ll}
1)\ &
{{[x_{3},\,[x_{2},\,[x_{1},\,x_{2}]]]}-
    {[[x_{1},\,x_{2}],\,[x_{2},\,x_{3}]]}}=0;\\
&
{{[[x_{1},\,x_{2}],\,[x_{2},\,x_{2}]]}=0};
  \;\;
  {{[[x_{2},\,x_{2}],\,[x_{2},\,x_{3}]]}=
  0},\\
&
 {{[[x_{3},\,[x_{1},\,x_{2}]],\,[x_{3},\,[x_{2},\,x_{2}]]]}= 0};\quad
 {{[[x_{3},\,[x_{1},\,x_{2}]],\,[[x_{1},\,x_{2}],\,[x_{1},\,x_{2}]]]}=0.}
\end{array}
$

$
\arraycolsep=0pt
\renewcommand{\arraystretch}{1.2}
\begin{array}{ll}
2)\ &  {{[[x_{1},\,x_{2}],\,[x_{2},
  \,[x_{2},\,x_{3}]]]}= 0};\quad
  {{[[x_{2},\,x_{3}],\,[x_{2},
  \,[x_{2},\,x_{3}]]]}= 0};\\
&  {{[[x_{3},\,[x_{1},\,x_{2}]],\,
      [x_{3},\,[x_{2},\,[x_{2},\,x_{3}]]]]}=0};\\
&  {{[[x_{3},\,[x_{1},\,x_{2}]],\,
      [[x_{1},\,x_{2}],\,[x_{3},
      \,[x_{1},\,x_{2}]]]]}=0.}
      \end{array}
$

\item[]

\hskip2.5ex\underline{$\fg(3, 6), \;\sdim =(36|40)$}

\item[]

$
\arraycolsep=0pt
\renewcommand{\arraystretch}{1.2}
\begin{array}{ll}
1)\ & {{[x_{4},\,[x_{3},\,[x_{2},\,x_{3}]]]}-
   {[[x_{2},\,x_{3}],\,[x_{3},\,x_{4}]]}}=0;\quad
 {{[[x_{2},\,x_{3}],\,[x_{3},\,x_{3}]]}=0};\\
 &
 {{[[x_{3},\,x_{3}],\,[x_{3},\,x_{4}]]}=0};\quad
 {{[[x_{4},\,[x_{2},\,x_{3}]],\,[x_{4},\,[x_{3},\,x_{3}]]]}=0};\\
&
 {{[[x_{4},\,[x_{2},\,x_{3}]],\,
   [[x_{2},\,x_{3}],\,[x_{2},\,x_{3}]]]}= 0};\\
& [[x_{3},\,[x_{2},\,x_{3}]],\,
   [[x_{1},\,x_{2}],\,[x_{3},\,x_{4}]]]{-}[[x_{4},\,[x_{2},\,x_{3}]],\,
   [[x_{1},\,x_{2}],\,[x_{3},\,x_{3}]]]{=}0.
\end{array}
$

$
\arraycolsep=0pt
\renewcommand{\arraystretch}{1.2}
\begin{array}{ll}
2)\ &
{{[x_{4},\,[x_{3},\,[x_{2},\,x_{3}]]]}-
   {[[x_{2},\,x_{3}],\,[x_{3},\,x_{4}]]}}=0;\\
&
  {{[[x_{1},\,x_{2}],\,[x_{2},\,x_{3}]]}=0};\quad
  {{[[x_{2},\,x_{3}],\,[x_{3},\,x_{3}]]}=0};\quad
  {{[[x_{3},\,x_{3}],\,[x_{3},\,x_{4}]]}=0};\\
&  {{[[x_{3},\,[x_{2},\,x_{3}]],\,
    [x_{4},\,[x_{2},\,x_{3}]]]}=0};\quad
  {{[[x_{4},\,[x_{2},\,x_{3}]],\,
    [x_{4},\,[x_{3},\,x_{3}]]]}=0};\\
&
{{[[[x_{1},\,x_{2}],\,[x_{3},\,x_{4}]],\,
    [[x_{1},\,x_{2}],\,[[x_{1},\,x_{2}],\,
    [x_{3},\,x_{3}]]]]}= 0.}
      \end{array}
$

$
\arraycolsep=0pt
\renewcommand{\arraystretch}{1.2}
\begin{array}{ll}
3)\ &  {{[[x_{2},\,x_{3}],\,
  [x_{3},\,[x_{3},\,x_{4}]]]}=0};\quad
  {{\left(\ad_{[x_{3},\,x_{4}]}\right)^2(x_{3})}=0}; \\
&  {{[[x_{4},\,[x_{2},\,x_{3}]],\,
      [x_{4},\,[x_{3},\,[x_{3},\,x_{4}]]]]}=0};\\
&  {{[[x_{4},\,[x_{2},\,x_{3}]],\,
      [[x_{2},\,x_{3}],\,
      [x_{4},\,[x_{2},\,x_{3}]]]]}=0};\\
&  {[[x_{4},\,[x_{2},\,x_{3}]],\,
      [[x_{3},\,x_{4}],\,[x_{3},\,
      [x_{1},\,x_{2}]]]]}-\\
&   {2\,
[[[x_{1},\,x_{2}],\,[x_{3},\,x_{4}]],\,
        [[x_{2},\,x_{3}],\,[x_{3},\,x_{4}]]]}=0.
        \end{array}
$

$
\arraycolsep=0pt
\renewcommand{\arraystretch}{1.2}
\begin{array}{ll}
5)\ & {{[[x_{1},\,x_{2}],\,[x_{2},\,x_{3}]]}=0};\quad
  {{[[x_{2},\,x_{3}],\,[x_{3},\,[x_{2},\,x_{3}]]]}=0};\\
& {{[[x_{3},\,x_{4}],\,[x_{3},\,[x_{2},\,x_{3}]]]}=0};\quad
  {{[[x_{4},\,[x_{2},\,x_{3}]],\,
    [[x_{2},\,x_{3}],\,[x_{2},\,x_{3}]]]}=0;}\\
&  {{[[x_{4},\,[x_{2},\,x_{3}]],\,
    [[x_{2},\,x_{3}],\,[x_{3},\,x_{4}]]]}=0};\\
& {{[[[x_{1},\,x_{2}],\,[x_{3},\,x_{4}]],\,
    [[x_{1},\,x_{2}],\,
    [x_{3},\,[x_{3},\,[x_{1},\,x_{2}]]]]]}=0.}
\end{array}
$

$ \arraycolsep=0pt
\renewcommand{\arraystretch}{1.2}
\begin{array}{ll}
4)\ & {{\left(\ad_{[x_{3},\,x_{4}]}\right)^2(x_{3})}\,{=}\,0};\;\;
 {{[[x_{1},\,x_{2}],\,[x_{2},\,x_{3}]]}\,{=}\,0};\;\;
 {{[[x_{2},\,x_{3}],\,[x_{3},\,[x_{3},\,x_{4}]]]}\,{=}\,0};\\
&{{[[x_{4},\,[x_{2},\,x_{3}]],\,[x_{4},\,
   [x_{3},\,[x_{3},\,x_{4}]]]]}\,{=}\,0};\\
& {{[[x_{4},\,[x_{2},\,x_{3}]],\,
   [[x_{2},\,x_{3}],\,[x_{3},\,x_{4}]]]}\,{=}\,0};\\
&
{{[[[x_{1},\,x_{2}],\,[x_{3},\,x_{4}]],\,
   [[x_{3},\,[x_{1},\,x_{2}]],\,
   [[x_{1},\,x_{2}],\,[x_{3},\,x_{4}]]]]}=0.}
\end{array}
$

$
\arraycolsep=0pt
\renewcommand{\arraystretch}{1.2}
\begin{array}{ll}
6)\ &
{{[[x_{1},\,x_{2}],\,[x_{2},\,x_{3}]]}=0};\quad
  {{[[x_{2},\,x_{3}],\,
  [[x_{2},\,x_{3}],\,[x_{3},\,x_{4}]]]}=0};\\
&
{{[[x_{3},\,x_{4}],\,[[x_{2},\,x_{3}],\,
  [x_{3},\,x_{4}]]]}\,{=}\,0};\\
& {{[[x_{4},\,[x_{2},\,x_{3}]],\,
      [[x_{2},\,x_{3}],\,[x_{4},\,
      [x_{2},\,x_{3}]]]]}\,{=}\,0};\\
& {{[[x_{4},\,[x_{2},\,x_{3}]],\,
      [[x_{3},\,x_{4}],\,[x_{4},\,
      [x_{2},\,x_{3}]]]]}=0};\\
&
{{[[[x_{1},\,x_{2}],\,[x_{3},\,x_{4}]],\,
      [[x_{3},\,[x_{1},\,x_{2}]],\,
      [[x_{1},\,x_{2}],\,[x_{3},\,x_{4}]]
       ]]}=0.}
       \end{array}
$

$
\arraycolsep=0pt
\renewcommand{\arraystretch}{1.2}
\begin{array}{ll}
7)\ & \left (\mathrm{ad}_{[x_2,x_4]}\right )^2([x_2,x_3])=0;\quad
\left (\mathrm{ad}_{[x_2,x_4]}\right )^2([x_3,x_4])=0;\\
&{{[x_{3},\,[x_{2},\,x_{4}]]}-{2\,
  [x_{4},\,[x_{2},\,x_{3}]]}=0};\\
& {{[[x_{4},\,[x_{1},\,x_{2}]],\,
      [[x_{3},\,[x_{1},\,x_{2}]],\,
      [x_{4},\,[x_{1},\,x_{2}]]]]}\,{=}\,0.}
       \end{array}
$

\end{arab}
}

\vskip .3cm

\underline{$\fg(3, 3), \;\sdim =(21|16)$} {\bf SuperLie} tells us
that $\sdim =(22|16)$: indeed the Cartan matrices below determine
a nontrivial central extension of a simple Lie superalgebra. The
relations below serve also the simple quotient algebras.

\footnotesize{
\begin{arab}

\item[]

$
\arraycolsep=0pt
\renewcommand{\arraystretch}{1.2}
\begin{array}{ll}
1)\ &
{{[[x_{3},\,x_{4}],\,[x_{3},\,[x_{2},\,x_{3}]]]}=0};\quad
  {{[[x_{3},\,[x_{1},\,x_{2}]],\,[x_{3},\,[x_{2},\,x_{3}]]]}=0};\\
&
{{[[[x_{1},\,x_{2}],\,[x_{3},\,x_{4}]],\,
    [[x_{2},\,x_{3}],\,[x_{3},\,x_{4}]]]}=0.}
\end{array}
$

$
\arraycolsep=0pt
\renewcommand{\arraystretch}{1.2}
\begin{array}{ll}
2)\ &
{{\left(\ad_{[x_{3},\,x_{4}]}\right)^2([x_{2},\,x_{3}])}=0};\quad
  {{[[x_{3},\,[x_{1},\,x_{2}]],\,
    [[x_{2},\,x_{3}],\,[x_{3},\,x_{4}]]]}=0};\\
&
{{[[[x_{1},\,x_{2}],\,[x_{3},\,x_{4}]],\,
    [[x_{3},\,x_{4}],\,[x_{4},\,[x_{2},\,x_{3}]]]]}=0.}
\end{array}
$

$
\arraycolsep=0pt
\renewcommand{\arraystretch}{1.2}
\begin{array}{ll}
3)\ &
{{[x_{3},\,[x_{2},\,x_{4}]]}-{[x_{4},\,[x_{2},\,x_{3}]]}=0};\quad
  {{[[x_{1},\,x_{2}],\,[x_{2},\,x_{3}]]}=0};\\
&
{{[[x_{1},\,x_{2}],\,[x_{2},\,x_{4}]]}=0};\quad
  {{[[x_{4},\,[x_{1},\,x_{2}]],\,
    [[x_{2},\,x_{4}],\,[x_{3},\,x_{4}]]]}=0.}
\end{array}
$

$
\arraycolsep=0pt
\renewcommand{\arraystretch}{1.2}
 \begin{array}{ll}
4)\ &
  {{[[x_{1},\,x_{2}],\,[x_{2},\,x_{3}]]}=0};\quad
  {{[[x_{1},\,x_{2}],\,[x_{2},\,x_{4}]]}=0};\\
& {{[[x_{2},\,x_{4}],\,[[x_{2},\,x_{3}],\,[x_{2},\,x_{4}]]]}\,{=}\,0};\\
& {{[[x_{4},\,[x_{1},\,x_{2}]],\,[[x_{2},\,x_{4}],\,
    [x_{4},\,[x_{2},\,x_{3}]]]]}\,{=}\,0.}
\end{array}
$

$
\arraycolsep=0pt
\renewcommand{\arraystretch}{1.2}
\begin{array}{ll}
5)\ &
  {{[[x_{4},\,[x_{1},\,x_{2}]],\,
      [[x_{3},\,[x_{1},\,x_{2}]],
      \,[x_{4},\,[x_{1},\,x_{2}]]]]}
    =0;}\\
 &   {{[[x_{4},\,[x_{1},\,x_{2}]],\,
      [[x_{4},\,[x_{1},\,x_{2}]],\,
      [x_{4},\,[x_{2},\,x_{3}]]]]}=0}.

\end{array}
$

$ \arraycolsep=0pt
\renewcommand{\arraystretch}{1.2}
\begin{array}{ll}
6)\ &
  {{[[x_{1},\,x_{2}],\,[x_{2},\,[x_{2},\,x_{3}]]]}=0};\quad
   {{[[x_{2},\,x_{4}],\,[x_{2},\,[x_{2},\,x_{3}]]]}=0;}\\
 & {{[[x_{2},\,x_{3}],\,[x_{4},\,[x_{1},\,x_{2}]]]}-
   {2\,
    [[x_{2},\,x_{4}],\,[x_{3},\,[x_{1},\,x_{2}]]]}=0.}
  \end{array}
$

$ \arraycolsep=0pt
\renewcommand{\arraystretch}{1.2}
\begin{array}{ll}
7)\ &
  {{[[x_{1},\,x_{2}],\,[x_{2},\,x_{4}]]}=0};\quad
  {{[[x_{4},\,[x_{1},\,x_{2}]],\,
      [[x_{1},\,x_{2}],\,[x_{3},\,[x_{1},\,x_{2}]]]]}=
      0};\\
  &{{[[x_{2},\,x_{3}],\,[x_{2},\,x_{4}]]}=0};
  \quad
  \left ( \mathrm{ad}_{[x_1,x_2]}\right )^2([x_2,x_3])=0; \quad
\end{array}
$

$ \arraycolsep=0pt
\renewcommand{\arraystretch}{1.2}
\begin{array}{ll}
8)\ &
  {{[x_{2},\,[x_{1},\,x_{3}]]}-{[x_{3},\,[x_{1},\,x_{2}]]}=0};\quad
  {{[[x_{1},\,x_{2}],\,[x_{2},\,x_{4}]]}=0};\\
  &
  {{[[x_{2},\,x_{3}],\,[x_{2},\,x_{4}]]}=0};\quad
  {{[[x_{4},\,[x_{1},\,x_{2}]],\,[[x_{1},\,x_{2}],\,
  [x_{1},\,x_{3}]]]}= 0.}
\end{array}
$

$ \arraycolsep=0pt
\renewcommand{\arraystretch}{1.2}
\begin{array}{ll}
9)\ &
  {{[[x_{1},\,x_{3}],\,[[x_{1},\,x_{3}],\,
    [x_{2},\,x_{3}]]]}=0};\quad
  {{[[x_{4},\,[x_{2},\,x_{3}]],\,
    [[x_{1},\,x_{3}],\,[x_{2},\,x_{3}]]]}=0};\\
 &{{[[[x_{1},\,x_{3}],\,[x_{2},\,x_{4}]],\,
    [[x_{1},\,x_{3}],\,[x_{2},\,[x_{1},\,x_{3}]]]]}=0.}
\end{array}
$

$ \arraycolsep=0pt
\renewcommand{\arraystretch}{1.2}
\begin{array}{ll}
10)\ &
   {{[[x_{1},\,x_{3}],\,[x_{3},\,[x_{2},\,x_{3}]]]}=0},
   {{[[x_{3},\,[x_{2},\,x_{3}]],\,
   [x_{4},\,[x_{2},\,x_{3}]]]}=0};\\
   & {{[[[x_{1},\,x_{3}],\,[x_{2},\,x_{3}]],\,
      [[x_{1},\,x_{3}],\,[x_{2},\,x_{4}]]]}=0}
      \end{array}
$
\item[]

\hskip2.5ex \underline{$\fg(4, 3), \;\sdim =(24|26)$}

\item[]

$
\arraycolsep=0pt
\renewcommand{\arraystretch}{1.2}
\begin{array}{ll}
1)\ &
{{[[x_{1},\,x_{2}],\,[x_{2},\,[x_{2},\,x_{3}]]]}=0};\\
&
{{[[x_{2},\,x_{3}],\,[x_{4},\,[x_{1},\,x_{2}]]]}-
  {2\,[[x_{2},\,x_{4}],\,[x_{3},\,[x_{1},\,x_{2}]]]}=0};\\
&
{{[[x_{2},\,x_{4}],\,[x_{2},\,[x_{2},\,x_{3}]]]}=0};\\
&  {{[[x_{4},\,[x_{1},\,x_{2}]],\,
   [[x_{1},\,x_{2}],\,[x_{3},\,[x_{1},\,x_{2}]]]]}=0.}
\end{array}
$

$
\arraycolsep=0pt
\renewcommand{\arraystretch}{1.2}
\begin{array}{ll}
2)\ &
{{[[x_{1},\,x_{2}],\,[x_{2},\,x_{3}]]}=0};\quad
  {{[[x_{2},\,x_{3}],\,[x_{2},\,x_{4}]]}=0};\\
&
{{[[x_{2},\,x_{4}],\,[[x_{1},\,x_{2}],\,
   [x_{2},\,x_{4}]]]}\,{=}\,0};\\
 & {{[[x_{4},\,[x_{2},\,x_{3}]],\,
   [[x_{2},\,x_{4}],\,[x_{4},\,
   [x_{1},\,x_{2}]]]]}\,{=}\,0};\\
& {{[[x_{4},\,[x_{2},\,x_{3}]],\,
   [[x_{3},\,[x_{1},\,x_{2}]],\,
   [x_{4},\,[x_{2},\,x_{3}]]]]}=0.}
\end{array}
$

$ \arraycolsep=0pt
\renewcommand{\arraystretch}{1.2}
\begin{array}{ll}
3)\ &
  {{[x_{3},\,[x_{2},\,x_{4}]]}-
    {[x_{4},\,[x_{2},\,x_{3}]]}=0};\quad
  {{[[x_{1},\,x_{2}],\,[x_{2},\,x_{3}]]}=0};\\
  &
  {{[[x_{1},\,x_{2}],\,[x_{2},\,x_{4}]]}=0};\quad
  {{[[x_{4},\,[x_{1},\,x_{2}]],\,
      [[x_{2},\,x_{4}],\,[x_{3},\,x_{4}]]]}=0};\\
 & {{[[x_{4},\,[x_{1},\,x_{2}]],\,
      [[x_{3},\,[x_{1},\,x_{2}]],\,
      [x_{4},\,[x_{1},\,x_{2}]]]]}=0}.
\end{array}
$

$ \arraycolsep=0pt
\renewcommand{\arraystretch}{1.2}
\begin{array}{ll}
4)\ &
  {{[x_{3},\,[x_{2},\,x_{4}]]}-
  {[x_{4},\,[x_{2},\,x_{3}]]}=0}; \quad
  {{[[x_{1},\,x_{2}],\,[x_{2},\,x_{3}]]}=0};\\
  &{{[[x_{1},\,x_{2}],\,[x_{2},\,x_{4}]]}=
  0}; \quad
  {{[[x_{4},\,[x_{1},\,x_{2}]],\,
      [[x_{2},\,x_{4}],\,[x_{3},\,x_{4}]]]}=0};\\
 & {{[[x_{4},\,[x_{1},\,x_{2}]],\,
      [[x_{3},\,[x_{1},\,x_{2}]],\,[x_{4},\,
      [x_{1},\,x_{2}]]]]}=0}.

\end{array}
$

$ \arraycolsep=0pt
\renewcommand{\arraystretch}{1.2}
\begin{array}{ll}
5)\ &
   {{[[x_{3},\,[x_{1},\,x_{2}]],\,
      [[x_{2},\,x_{3}],\,[x_{3},\,x_{4}]]]}=0};\\
& {{[[[x_{1},\,x_{2}],\,[x_{3},\,x_{4}]],\,
      [[x_{3},\,x_{4}],\,[x_{4},\,[x_{2},\,x_{3}]]]]}=0};\\
 & {{[[[x_{1},\,x_{2}],\,[x_{3},\,x_{4}]],\,
      [[x_{3},\,[x_{1},\,x_{2}]],\,
       [[x_{1},\,x_{2}],\,[x_{3},\,x_{4}]]]]}=0.}

\end{array}
$

$ \arraycolsep=0pt
\renewcommand{\arraystretch}{1.2}
\begin{array}{ll}
6)\ & %
  {{[[x_{1},\,x_{2}],\,[x_{2},\,x_{3}]]}= 0};\quad
   \left (\mathrm{ad}_{[x_2,x_3]} \right )^2([x_3,x_4])=0; \\
   &{{[[x_{3},\,x_{4}],\,
      [[x_{2},\,x_{3}],\,[x_{3},\,x_{4}]]]}=0};\\
 &  {{[[x_{4},\,[x_{2},\,x_{3}]],\,
      [[x_{3},\,x_{4}],\,[x_{4},\,[x_{2},\,x_{3}]]]]}=0;}\\
      &
  {{[[x_{4},\,[x_{2},\,x_{3}]],\,
      [[x_{2},\,x_{3}],\,[x_{4},\,[x_{2},\,x_{3}]]]]}=0.}

\end{array}
$

$ \arraycolsep=0pt
\renewcommand{\arraystretch}{1.2}
\begin{array}{ll}
7)\ &
  {{[[x_{3},\,x_{4}],\,[x_{3},\,[x_{2},\,x_{3}]]]}=0};\quad
  {{[[x_{3},\,[x_{1},\,x_{2}]],\,
  [x_{3},\,[x_{2},\,x_{3}]]]}=0};\\
  &
  {{[[[x_{1},\,x_{2}],\,[x_{3},\,x_{4}]],\,
      [[x_{2},\,x_{3}],\,[x_{3},\,x_{4}]]]}=0};\\
      &
  {{[[[x_{1},\,x_{2}],\,[x_{3},\,x_{4}]],\,
    [[x_{1},\,x_{2}],\,[x_{3},\,[x_{3},\,[x_{1},\,x_{2}]]]]]}= 0.}

\end{array}
$

$ \arraycolsep=0pt
\renewcommand{\arraystretch}{1.2}
\begin{array}{ll}
8)\ &
  {{[[x_{1},\,x_{2}],\,[x_{2},\,x_{3}]]}=0};\quad
  {{[[x_{2},\,x_{3}],\,[x_{3},\,[x_{3},\,x_{4}]]]}=0;}\\
  &
  \left (\mathrm{ad}_{[x_3,x_4]}\right )^2 (x_3)=0; \quad
  {{[[x_{4},\,[x_{2},\,x_{3}]],\,
      [x_{4},\,[x_{3},\,[x_{3},\,x_{4}]]]]}=0};\\
      &
  {{[[x_{4},\,[x_{2},\,x_{3}]],\,
      [[x_{2},\,x_{3}],\,[x_{3},\,x_{4}]]]}=0.}

\end{array}
$

$ \arraycolsep=0pt
\renewcommand{\arraystretch}{1.2}
\begin{array}{ll}
9)\ & 
  {{[[x_{1},\,x_{2}],\,[x_{2},\,x_{3}]]}=0};\quad
  \left ( \mathrm{ad}_{[x_2,x_3]}\right )^2(x_3)=0;\\
  &
  {{[[x_{3},\,x_{4}],\,[x_{3},\,[x_{2},\,x_{3}]]]}=0};\quad
  {{[[x_{4},\,[x_{2},\,x_{3}]],\,
      [[x_{2},\,x_{3}],\,[x_{2},\,x_{3}]]]}=0};\\
      &
  {{[[x_{4},\,[x_{2},\,x_{3}]],\,
      [[x_{2},\,x_{3}],\,[x_{3},\,x_{4}]]]}=0.}
\end{array}
$

$
\arraycolsep=0pt
\renewcommand{\arraystretch}{1.2}
\begin{array}{ll}
10)\ &
  {{[x_{3},\,[x_{3},\,x_{3}]]}=0};\quad
  {{[x_{4},\,[x_{3},\,[x_{2},\,x_{3}]]]}-
  {[[x_{2},\,x_{3}],\,[x_{3},\,x_{4}]]}=0};\\
  &
  {{[[x_{1},\,x_{2}],\,[x_{2},\,x_{3}]]}=
  0};\quad
  {{[[x_{2},\,x_{3}],\,[x_{3},\,x_{3}]]}=0};\\
  &
  {{[[x_{3},\,x_{3}],\,[x_{3},\,x_{4}]]}=0};\quad
  {{[[x_{3},\,[x_{2},\,x_{3}]],\,[x_{4},\,
  [x_{2},\,x_{3}]]]}=0};\\
  &{{[[x_{4},\,[x_{2},\,x_{3}]],\,[x_{4},\,[x_{3},\,x_{3}]]]}
     =0.}

\end{array}
$

\item[]

\hskip2.5ex \underline{$\fg(2, 6), \; \sdim =(35|20)$}\vspace{-2mm}

\item[]

$
\arraycolsep=0pt
\renewcommand{\arraystretch}{1.2}
\begin{array}{ll}
1)\ &
{{[[x_{2},\,x_{3}],\,[x_{3},\,x_{5}]]}=0};\quad
  {{[[x_{3},\,x_{4}],\,[x_{3},\,x_{5}]]}=0};\\
&  {{[[x_{3},\,[x_{1},\,x_{2}]],\,
    [[x_{2},\,x_{3}],\,[x_{3},\,x_{4}]]]}=0.}
\end{array}
$
\vskip -1mm

$
\arraycolsep=0pt
\renewcommand{\arraystretch}{1.2}
\begin{array}{ll}
2)\ &
  {{[x_{3},\,[x_{2},\,x_{4}]]}-{[x_{4},\,[x_{2},\,x_{3}]]}=0};\quad
  {{[[x_{1},\,x_{2}],\,[x_{2},\,x_{3}]]}=0};\\
  &{{[[x_{1},\,x_{2}],\,[x_{2},\,x_{4}]]}=0};\quad
  {{[[x_{2},\,x_{3}],\,[x_{3},\,x_{5}]]}=0};\quad
  {{[[x_{3},\,x_{4}],\,[x_{3},\,x_{5}]]}=0.}
\end{array}
$
\vskip -1mm

$ \arraycolsep=0pt
\renewcommand{\arraystretch}{1.2}
\begin{array}{ll}
3)\
  {{[[[x_{1},\,x_{2}],\,[x_{3},\,x_{5}]],\,
      [[x_{5},\,[x_{2},\,x_{3}]],\,[x_{5},\,
      [x_{3},\,x_{4}]]]]}=0.}
     \end{array}
     $
\vskip -1mm

$ \arraycolsep=0pt
\renewcommand{\arraystretch}{1.2}
\begin{array}{ll}
4)\ &
  {{[[x_{1},\,x_{2}],\,[x_{2},\,x_{3}]]}=0};\quad
  {{[[x_{1},\,x_{2}],\,[x_{2},\,x_{4}]]}=0};\\
&  {{[[x_{5},\,[x_{2},\,x_{3}]],\,
  [[x_{2},\,x_{3}],\,[x_{2},\,x_{4}]]
      ]}=0.}
\end{array}
$
\vskip -1mm

$ \arraycolsep=0pt
\renewcommand{\arraystretch}{1.2}
\begin{array}{ll}
5)\ &
  {{[[x_{4},\,[x_{1},\,x_{2}]],\,
  [[x_{2},\,x_{4}],\,[x_{3},\,x_{4}]]
      ]}=0};\\
      &{{[[x_{5},\,[x_{3},\,x_{4}]],\,
      [[x_{2},\,x_{4}],\,[x_{3},\,x_{4}]]]}=0.}
      \end{array}
      $
\vskip -1mm

$ \arraycolsep=0pt
\renewcommand{\arraystretch}{1.2}
\begin{array}{ll}
6)\ &
  {{[[[x_{1},\,x_{2}],\,[x_{3},\,x_{5}]],\,
      [[x_{3},\,[x_{1},\,x_{2}]],\,[x_{4},\,
      [x_{1},\,x_{2}]]]]}=0.}
     \end{array}
     $
\vskip -1mm


\item[]

\hskip2.5ex \underline{$\fg(8, 3), \;\sdim=(55|50)$}\vspace{-2mm}

\enlargethispage{2mm}

\item[]

$
\arraycolsep=0pt
\renewcommand{\arraystretch}{1.2}
\begin{array}{ll}
1)\ &
{{[x_{3},\,[x_{3},\,[x_{2},\,x_{3}]]]}=0};\quad
  {{[[x_{3},\,x_{4}],\,[x_{3},\,[x_{2},\,x_{3}]]]}=0};\\
& {{[[x_{3},\,[x_{1},\,x_{2}]],\,[x_{3},\,[x_{2},\,x_{3}]]]}=0};\\
&
{{[[[x_{1},\,x_{2}],\,[x_{3},\,x_{4}]],\,
    [[x_{3},\,x_{4}],\,[x_{4},\,[x_{2},\,x_{3}]]]]}=0};\\
&
{[[[x_{1},\,x_{2}],\,[x_{3},\,x_{4}]],\,
    [[x_{4},\,x_{5}],\,[x_{3},\,[x_{2},\,x_{3}]]]]}-\\
&
{[[[x_{2},\,x_{3}],\,[x_{4},\,x_{5}]],\,
    [[x_{3},\,x_{4}],\,[x_{3},\,[x_{1},\,x_{2}]]]]}=0.
\end{array}
$

$
\arraycolsep=0pt
\renewcommand{\arraystretch}{1.2}
\begin{array}{ll}
2)\ &
{{[[x_{3},\,x_{4}],\,[x_{4},\,x_{5}]]}=0};\quad
  {{[[x_{3},\,x_{4}],\,
    [x_{3},\,[x_{2},\,x_{3}]]]}=0};\\
&  {{[[x_{3},\,[x_{1},\,x_{2}]],\,
    [x_{3},\,[x_{2},\,x_{3}]]]}\,{=}\,0};
    \\
&
{{[[[x_{1},\,x_{2}],\,[x_{3},\,x_{4}]],\,
    [[x_{2},\,x_{3}],\,[x_{3},\,x_{4}]]]}\,{=}\,0};\\
&
{{[[[x_{4},\,x_{5}],\,[x_{3},\,[x_{1},\,x_{2}]]],\,
      [[x_{5},\,[x_{3},\,x_{4}]],\,
       [[x_{2},\,x_{3}],\,[x_{4},\,x_{5}]]]]}=0.}
\end{array}
$

$ \arraycolsep=0pt
\renewcommand{\arraystretch}{1.2}
\begin{array}{ll}
3)\ &
{{[[x_{3},\,x_{4}],\,[x_{4},\,x_{5}]]}=0};\quad
  {{[[x_{3},\,x_{4}],\,[[x_{2},\,x_{3}],\,[x_{3},\,x_{4}]]]}
     =0};\\
     &{{[[x_{3},\,[x_{1},\,x_{2}]],\,
      [[x_{2},\,x_{3}],\,[x_{3},\,x_{4}]]]}=0};\\
& {{[[[x_{1},\,x_{2}],\,[x_{3},\,x_{4}]],\,
      [[x_{3},\,x_{4}],\,[x_{4},\,
      [x_{2},\,x_{3}]]]]}=0};\\
 & {{[[[x_{4},\,x_{5}],\,[x_{3},\,[x_{1},\,x_{2}]]],\,
      [[x_{5},\,[x_{3},\,x_{4}]],\,
       [[x_{2},\,x_{3}],\,[x_{4},\,x_{5}]]]]}=0.}
\end{array}
$

$ \arraycolsep=0pt
\renewcommand{\arraystretch}{1.2}
\begin{array}{ll}
4)\ &
  {{[x_{3},\,[x_{2},\,x_{4}]]}-{[x_{4},\,[x_{2},\,x_{3}]]}=0};\quad
  {{[[x_{1},\,x_{2}],\,[x_{2},\,x_{3}]]}=0};\\
 & {{[[x_{1},\,x_{2}],\,[x_{2},\,x_{4}]]}=0};\quad
  {{[[x_{4},\,[x_{1},\,x_{2}]],\,
      [[x_{2},\,x_{4}],\,[x_{3},\,x_{4}]]]}=0};\\
      &
  {{[[[x_{1},\,x_{2}],\,[x_{4},\,x_{5}]],\,
      [[x_{5},\,[x_{2},\,x_{4}]],\,
      [x_{5},\,[x_{3},\,x_{4}]]]
      ]}=0.}
\end{array}
$

$ \arraycolsep=0pt
\renewcommand{\arraystretch}{1.2}
\begin{array}{ll}
5)\ &
  {{[[x_{1},\,x_{2}],\,[x_{2},\,x_{3}]]}=0};\quad
  {{[[x_{1},\,x_{2}],\,[x_{2},\,x_{4}]]}=0};\\
  &
  {{[[x_{2},\,x_{4}],\,[x_{4},\,x_{5}]]}=0};\quad
 \left (\mathrm{ad}_{[x_2,x_4]} \right )^2([x_2,x_3])=0;\\
     & {{[[x_{4},\,[x_{1},\,x_{2}]],\,
      [[x_{2},\,x_{4}],\,[x_{4},\,[x_{2},\,x_{3}]]]]}=0};\\
& {{[[[x_{1},\,x_{2}],\,[x_{4},\,x_{5}]],\,
      [[x_{5},\,[x_{2},\,x_{4}]],\,
       [[x_{2},\,x_{3}],\,[x_{4},\,x_{5}]]]]}=0.}
\end{array}
$

$ \arraycolsep=0pt
\renewcommand{\arraystretch}{1.2}
\begin{array}{ll}
6)\ &
  {{[[x_{2},\,x_{4}],\,[x_{4},\,x_{5}]]}=0};\quad
  {{[[x_{4},\,[x_{1},\,x_{2}]],\,
      [[x_{3},\,[x_{1},\,x_{2}]],\,[x_{4},\,[x_{1},\,x_{2}]]]
      ]}=0};\\
      & {{[[x_{4},\,[x_{1},\,x_{2}]],\,
      [[x_{4},\,[x_{1},\,x_{2}]],\,[x_{4},\,[x_{2},\,x_{3}]]]
      ]}=0};\\
      &{{[[[x_{1},\,x_{2}],\,[x_{4},\,x_{5}]],\,
      [[[x_{1},\,x_{2}],\,[x_{4},\,x_{5}]],\,
       [[x_{2},\,x_{3}],\,[x_{4},\,x_{5}]]]]}=0.}
       \end{array}
       $

       $
        \arraycolsep=0pt
\renewcommand{\arraystretch}{1.2}
\begin{array}{ll}
7)\ &
  {{[[x_{2},\,x_{4}],\,[x_{4},\,x_{5}]]}=0};\quad
  {{[[x_{1},\,x_{2}],\,[x_{2},\,[x_{2},\,x_{3}]]]}=0};\\
  &
  {{[[x_{2},\,x_{3}],\,[x_{4},\,[x_{1},\,x_{2}]]]}-
    {2\,
    [[x_{2},\,x_{4}],\,[x_{3},\,[x_{1},\,x_{2}]]]}=0;}\\
    &
  {{[[x_{2},\,x_{4}],\,[x_{2},\,[x_{2},\,x_{3}]]]}=0;}\\
 & {{[[[x_{1},\,x_{2}],\,[x_{4},\,x_{5}]],\,
      [[x_{5},\,[x_{2},\,x_{4}]],\,
       [[x_{2},\,x_{3}],\,[x_{4},\,x_{5}]]]]}=0.}
       \end{array}$

       $
       \arraycolsep=0pt
\renewcommand{\arraystretch}{1.2}
\begin{array}{ll}
8)\ &
  {{[[x_{1},\,x_{2}],\,[x_{2},\,x_{4}]]}=0};\quad
  {{[[x_{2},\,x_{3}],\,[x_{2},\,x_{4}]]}=0;}\\
&
{{[[x_{2},\,x_{4}],\,[x_{4},\,x_{5}]]}=0};\quad
  {{[[x_{1},\,x_{2}],\,[[x_{1},\,x_{2}],\,
  [x_{2},\,x_{3}]]]}=0;}\\
  & {{[[x_{4},\,[x_{1},\,x_{2}]],\,
      [[x_{1},\,x_{2}],\,[x_{3},\,[x_{1},\,x_{2}]]]]}=0;}\\
& {{[[[x_{1},\,x_{2}],\,[x_{4},\,x_{5}]],\,
      [[[x_{1},\,x_{2}],\,[x_{4},\,x_{5}]],\,
       [[x_{2},\,x_{3}],\,[x_{4},\,x_{5}]]]]}=0.}
      \end{array} $

      $
 \arraycolsep=0pt
\renewcommand{\arraystretch}{1.2}
\begin{array}{ll}
9)\ &
  {{[[x_{1},\,x_{2}],\,[x_{2},\,[x_{2},\,x_{3}]]]}=0};\\
 &{{[[x_{2},\,x_{3}],\,[x_{4},\,[x_{1},\,x_{2}]]]}-
    {2\,[[x_{2},\,x_{4}],\,[x_{3},\,[x_{1},\,x_{2}]]]}=0};\\
&
{{[[x_{2},\,x_{4}],\,[x_{2},\,[x_{2},\,x_{3}]]]}=0};\quad
  {{[[x_{4},\,[x_{1},\,x_{2}]],\,
      [[x_{2},\,x_{4}],\,[x_{4},\,[x_{2},\,x_{3}]]]]}=0.}
      \end{array}
      $

      $
\arraycolsep=0pt
\renewcommand{\arraystretch}{1.2}
\begin{array}{ll}
10)\ &
  {{[x_{2},\,[x_{1},\,x_{3}]]}-{[x_{3},\,[x_{1},\,x_{2}]]}=0};\quad
  {{[[x_{1},\,x_{2}],\,[x_{2},\,x_{4}]]}=0};\\
&  {{[[x_{4},\,[x_{1},\,x_{2}]],\,
      [[x_{1},\,x_{2}],\,[x_{1},\,x_{3}]]]}=0;}\quad
  {{[[x_{2},\,x_{3}],\,[x_{2},\,x_{4}]]}=0;}\\
& {{[[[x_{1},\,x_{2}],\,[x_{4},\,x_{5}]],\,
      [[[x_{1},\,x_{2}],\,[x_{4},\,x_{5}]],\,
       [[x_{2},\,x_{3}],\,[x_{4},\,x_{5}]]]]}=0.}
       \end{array}$

$  \arraycolsep=0pt
\renewcommand{\arraystretch}{1.2}
\begin{array}{ll}
11)\ &
{{[[x_{1},\,x_{2}],\,[x_{2},\,x_{4}]]}=0};\quad
  {{[[x_{2},\,x_{3}],\,[x_{2},\,x_{4}]]}=0};\\
 & \left ( \mathrm{ad}_{[x_1,x_2]}\right )^2([x_2,x_3])=0; \quad
  {{[[x_{4},\,[x_{1},\,x_{2}]],\,
      [[x_{1},\,x_{2}],\,[x_{3},\,[x_{1},\,x_{2}]]]]}=0};\\
&  {{[[x_{4},\,[x_{1},\,x_{2}]],\,
      [[x_{4},\,[x_{1},\,x_{2}]],\,
      [x_{4},\,[x_{2},\,x_{3}]]]]}=0.}
     \end{array}
     $

$ \arraycolsep=0pt
\renewcommand{\arraystretch}{1.2}
\begin{array}{ll}
12)\ &
  {{[[x_{1},\,x_{3}],\,[[x_{1},\,x_{3}],\,
  [x_{2},\,x_{3}]]]}=0};\;\;
  {{[[x_{4},\,[x_{2},\,x_{3}]],\,
      [[x_{1},\,x_{3}],\,[x_{2},\,x_{3}]]]}=0};\\
& {{[[[x_{1},\,x_{3}],\,[x_{2},\,x_{4}]],\,
      [[x_{1},\,x_{3}],\,[x_{2},\,[x_{1},\,x_{3}]]]]}=0};\\
&
{{[[[x_{4},\,x_{5}],\,[x_{2},\,[x_{1},\,x_{3}]]],\,
      [[[x_{2},\,x_{3}],\,[x_{4},\,x_{5}]],\,
       [[x_{4},\,x_{5}],\,[x_{2},\,[x_{1},\,x_{3}]]]]]}\,{=}\,0.}
   \end{array}
   $

$\arraycolsep=0pt
\renewcommand{\arraystretch}{1.2}
\begin{array}{ll}
13)\ &
  {{[x_{2},\,[x_{1},\,x_{3}]]}-{[x_{3},\,[x_{1},\,x_{2}]]}=0};\quad
  {{[[x_{1},\,x_{2}],\,[x_{2},\,x_{4}]]}=0};\\
 &{{[[x_{2},\,x_{3}],\,[x_{2},\,x_{4}]]}=0};\quad
 {{[[x_{2},\,x_{4}],\,[x_{4},\,x_{5}]]}=0};\\
 & {{[[x_{4},\,[x_{1},\,x_{2}]],\,
      [[x_{1},\,x_{2}],\,[x_{1},\,x_{3}]]]}=0};\\
&  {{[[x_{4},\,[x_{1},\,x_{2}]],\,
      [[x_{4},\,[x_{1},\,x_{2}]],\,
      [x_{4},\,[x_{2},\,x_{3}]]]]}=0.}
     \end{array}
     $

     $
\arraycolsep=0pt
\renewcommand{\arraystretch}{1.2}
\begin{array}{ll}
14)\ &
  {{[[x_{1},\,x_{3}],\,[x_{3},\,[x_{2},\,x_{3}]]]}=0};\quad
  {{[[x_{3},\,[x_{2},\,x_{3}]],\,
  [x_{4},\,[x_{2},\,x_{3}]]]}=0};\\
  &{{[[[x_{1},\,x_{3}],\,[x_{2},\,x_{3}]],\,
      [[x_{1},\,x_{3}],\,[x_{2},\,x_{4}]]]}=0};\\
&
{{[[[x_{4},\,x_{5}],\,[x_{2},\,[x_{1},\,x_{3}]]],\,
      [[x_{5},\,[x_{2},\,x_{4}]],\,
       [[x_{4},\,x_{5}],\,[x_{3},\,
       [x_{2},\,x_{3}]]]]]}=0.}
       \end{array}
       $

       $
       \arraycolsep=0pt
\renewcommand{\arraystretch}{1.2}
\begin{array}{ll}
15)\ &
  {{[[x_{2},\,x_{4}],\,[x_{4},\,x_{5}]]}=0};\quad
  {{[[x_{1},\,x_{3}],\,[[x_{1},\,x_{3}],\,
  [x_{2},\,x_{3}]]]}=0};\\
  &{{[[x_{4},\,[x_{2},\,x_{3}]],\,
      [[x_{1},\,x_{3}],\,[x_{2},\,x_{3}]]]}=0};\\
      &
  {{[[[x_{1},\,x_{3}],\,[x_{2},\,x_{4}]],\,
      [[x_{1},\,x_{3}],\,[x_{2},\,[x_{1},\,x_{3}]]]]}=0};\\
& {{[[[x_{1},\,x_{3}],\,[x_{2},\,x_{4}]],\,
      [[x_{4},\,[x_{2},\,x_{3}]],\,
       [[x_{1},\,x_{3}],\,[x_{2},\,x_{4}]]]]}=0.}
       \end{array}
       $

       $
        \arraycolsep=0pt
\renewcommand{\arraystretch}{1.2}
\begin{array}{ll}
16)\ &
     {{[x_{2},\,[x_{1},\,x_{3}]]}-
     {2\, [x_{3},\,[x_{1},\,x_{2}]]}=0}; \quad
     {{[[x_{2},\,x_{4}],\,[x_{4},\,x_{5}]]}=0};\\
 & \left (\mathrm{ad}_{[x_1,x_2]} \right )^2([x_1,x_3])=0;\quad
   \left (\mathrm{ad}_{[x_1,x_2]} \right )^2([x_2,x_3])=0.
\end{array}
     $

     $
 \arraycolsep=0pt
\renewcommand{\arraystretch}{1.2}
\begin{array}{ll}
17)\ &
  {{[[x_{2},\,x_{4}],\,[x_{4},\,x_{5}]]}=0};\quad
  {{[[x_{1},\,x_{3}],\,[x_{3},\,[x_{2},\,x_{3}]]]}=0};\\
 & {{[[x_{3},\,[x_{2},\,x_{3}]],\,
  [x_{4},\,[x_{2},\,x_{3}]]]}=0};\\
  & {{[[[x_{1},\,x_{3}],\,[x_{2},\,x_{3}]],\,
      [[x_{1},\,x_{3}],\,[x_{2},\,x_{4}]]]}=0};\\
& {{[[[x_{1},\,x_{3}],\,[x_{2},\,x_{4}]],\,
      [[x_{2},\,x_{4}],\,[x_{4},\,[x_{3},\,
      [x_{2},\,x_{3}]]]]]}=0.}
     \end{array}
     $

     $
     \arraycolsep=0pt
\renewcommand{\arraystretch}{1.2}
\begin{array}{ll}
18)\ &
  {{[[x_{2},\,x_{3}],\,[x_{2},\,x_{4}]]}=0};\quad
  {{[[x_{2},\,x_{4}],\,[x_{4},\,x_{5}]]}=0};\\
 & \left (\mathrm{ad}_{[x_1,x_3]} \right )^2([x_2,x_3])=0;\quad
   \left (\mathrm{ad}_{[x_2,x_3]}\right )^2([x_1,x_3])=0;\\
 & {{[[x_{2},\,[x_{1},\,x_{3}]],\,
     [[x_{1},\,x_{3}],\,[x_{2},\,[x_{1},\,x_{3}]]]]}=0};\\
&  {{[[x_{2},\,[x_{1},\,x_{3}]],\,
     [[x_{2},\,x_{3}],\,[x_{2},\,[x_{1},\,x_{3}]]]]}=0.}
      \end{array}
      $

      $
       \arraycolsep=0pt
\renewcommand{\arraystretch}{1.2}
\begin{array}{ll}
19)\ &
  {{[[x_{2},\,x_{3}],\,[x_{2},\,x_{4}]]}=0};\quad
  {{[[x_{2},\,x_{4}],\,[x_{4},\,x_{5}]]}=0};\\
& {{[[x_{1},\,x_{3}],\,[x_{3},\,[x_{2},\,x_{3}]]]}=0};\quad
  \left ( \mathrm{ad}_{[x_2,x_3]}\right )^2(x_3)=0; \\
& {{[[x_{2},\,[x_{1},\,x_{3}]],\,
      [[x_{1},\,x_{3}],\,[x_{2},\,x_{3}]]]}=0};\\
& {{[[x_{2},\,[x_{1},\,x_{3}]],\,
      [[x_{2},\,x_{3}],\,[x_{2},\,x_{3}]]]}=0.}
      \end{array}
      $

      $
       \arraycolsep=0pt
\renewcommand{\arraystretch}{1.2}
\begin{array}{ll}
20)\ &
  {{[[x_{2},\,x_{3}],\,[x_{2},\,x_{4}]]}=0};\quad
  \left (\mathrm{ad}_{[ x_1,x_3 ]} \right )^2(x_3)=0; \\
& {{[[x_{2},\,x_{3}],\,[x_{3},\,[x_{1},\,x_{3}]]]}=0};\quad
  {{[[x_{2},\,[x_{1},\,x_{3}]],\,[[x_{1},\,x_{3}],\,[x_{1},\,x_{3}]]]}=0};\\
& {{[[x_{2},\,[x_{1},\,x_{3}]],\,[[x_{1},\,x_{3}],\,[x_{2},\,x_{3}]]]}=0.}
\end{array}
      $

$
 \arraycolsep=0pt
\renewcommand{\arraystretch}{1.2}
\begin{array}{ll}
21)\ & {{[x_{3},\,[x_{3},\,x_{3}]]}=0};\quad
  {{[x_{3},\,[x_{2},\,[x_{1},\,x_{3}]]]}-
    {[[x_{1},\,x_{3}],\,[x_{2},\,x_{3}]]}=0};\\&
  {{[[x_{1},\,x_{3}],\,[x_{3},\,x_{3}]]}=0};\;\;
  {{[[x_{2},\,x_{3}],\,[x_{2},\,x_{4}]]}=0};\;\;
  {{[[x_{2},\,x_{3}],\,[x_{3},\,x_{3}]]}=0};\\&
  {{[[x_{2},\,[x_{1},\,x_{3}]],\,
  [x_{3},\,[x_{1},\,x_{3}]]]}=0};\quad
  {{[[x_{2},\,[x_{1},\,x_{3}]],\,[x_{3},\,[x_{2},\,x_{3}]]]}
=0.}
     \end{array}
     $

\newpage

\item[]

\hskip2.5ex \underline{$\fg(4, 6), \;\sdim =(66|32)$}

\item[]

$
\arraycolsep=0pt
\renewcommand{\arraystretch}{1.2}
\begin{array}{ll}
1)\ &
{{[[x_{3},\,x_{4}],\,[x_{4},\,x_{6}]]}=0};
\quad
  {{[[x_{4},\,x_{5}],\,[x_{4},\,x_{6}]]}=0};\\
&  {{[[x_{4},\,[x_{2},\,x_{3}]],\,
    [[x_{3},\,x_{4}],\,[x_{4},\,x_{5}]]]}=0.}
\end{array}
$

$
\arraycolsep=0pt
\renewcommand{\arraystretch}{1.2}
\begin{array}{ll}
2)\ &
{{[x_{4},\,[x_{3},\,x_{5}]]}-{[x_{5},\,[x_{3},\,x_{4}]]}=0};\quad
  {{[[x_{2},\,x_{3}],\,[x_{3},\,x_{4}]]}=0};\\
&
{{[[x_{2},\,x_{3}],\,[x_{3},\,x_{5}]]}=0};\quad
  {{[[x_{3},\,x_{4}],\,[x_{4},\,x_{6}]]}=0};\quad
{{[[x_{4},\,x_{5}],\,[x_{4},\,x_{6}]]}=0}.
\end{array}
$

$
\arraycolsep=0pt
\renewcommand{\arraystretch}{1.2}
\begin{array}{ll}
3)\ &
{{[[[x_{2},\,x_{3}],\,[x_{4},\,x_{6}]],\,
    [[x_{6},\,[x_{3},\,x_{4}]],\,
    [x_{6},\,[x_{4},\,x_{5}]]]]}=0}.
\end{array}
$

$\arraycolsep=0pt
\renewcommand{\arraystretch}{1.2}
\begin{array}{ll}
4)\ &
  {{[[x_{1},\,x_{2}],\,[x_{2},\,x_{3}]]}=0};\quad
  {{[[x_{2},\,x_{3}],\,[x_{3},\,x_{4}]]}=0};\quad
  {{[[x_{2},\,x_{3}],\,[x_{3},\,x_{5}]]}=0};\\&
  {{[[x_{6},\,[x_{3},\,x_{4}]],\,
  [[x_{3},\,x_{4}],\,[x_{3},\,x_{5}]]
      ]}=0.}
      \end{array}
      $

$\arraycolsep=0pt
\renewcommand{\arraystretch}{1.2}
\begin{array}{ll}
5)\ &
  {{[[x_{5},\,[x_{2},\,x_{3}]],\,
  [[x_{3},\,x_{5}],\,[x_{4},\,x_{5}]]]}=0};\\
 & {{[[x_{6},\,[x_{4},\,x_{5}]],\,[[x_{3},\,x_{5}],\,[x_{4},\,x_{5}]]]}=0.}
\end{array}
      $

$\arraycolsep=0pt
\renewcommand{\arraystretch}{1.2}
\begin{array}{ll}
6)\ &
  {{[[x_{1},\,x_{2}],\,[x_{2},\,x_{3}]]}=0};\\
&  {{[[[x_{2},\,x_{3}],\,[x_{4},\,x_{6}]],\,
     [[x_{4},\,[x_{2},\,x_{3}]],\,[x_{5},\,
     [x_{2},\,x_{3}]]]]}=0.}
\end{array}
     $

     $
\arraycolsep=0pt
\renewcommand{\arraystretch}{1.2}
\begin{array}{ll}
7)\ &
  {{[[[x_{4},\,x_{6}],\,[x_{3},\,[x_{1},\,x_{2}]]],\,
      [[[x_{1},\,x_{2}],\,[x_{3},\,x_{4}]],\,
       [[x_{1},\,x_{2}],\,[x_{3},\,x_{5}]]]]}=0.}
     \end{array}  $

\item[]

\hskip2.5ex \underline{$\fg(6, 6),\;\sdim =(78|64)$}

\item[]

$
\arraycolsep=0pt
\renewcommand{\arraystretch}{1.2}
\begin{array}{ll}
1)\ &
{{[[x_{1},\,x_{2}],\,[x_{2},\,x_{3}]]}=0};\quad
  {{[[x_{3},\,x_{4}],\,[x_{4},\,x_{6}]]}=0};\quad
  {{[[x_{4},\,x_{5}],\,[x_{4},\,x_{6}]]}=0};\\
&  {{[[x_{4},\,[x_{2},\,x_{3}]],\,
    [[x_{3},\,x_{4}],\,[x_{4},\,x_{5}]]]}=0};\\
&
{{[[[x_{2},\,x_{3}],\,[x_{4},\,x_{6}]],\,
    [[x_{4},\,[x_{2},\,x_{3}]],\,[[x_{2},\,x_{3}],\,
    [x_{4},\,x_{5}]]]]}=0}.
\end{array}
$

$
\arraycolsep=0pt
\renewcommand{\arraystretch}{1.2}
\begin{array}{ll}
2)\ &
{{[[x_{3},\,x_{4}],\,[x_{4},\,x_{6}]]}=0};\quad
  {{[[x_{4},\,x_{5}],\,[x_{4},\,x_{6}]]}=0};\\
&  {{[[x_{4},\,[x_{2},\,x_{3}]],\,
    [[x_{3},\,x_{4}],\,[x_{4},\,x_{5}]]]}=0};\\
&
{{[[[x_{4},\,x_{6}],\,[x_{3},\,[x_{1},\,x_{2}]]],\,
    [[[x_{1},\,x_{2}],\,[x_{3},\,x_{4}]],\,
    [[x_{4},\,x_{5}],\,[x_{3},\,
    [x_{1},\,x_{2}]]]]]}=0}.
\end{array}
$

$
\arraycolsep=0pt
\renewcommand{\arraystretch}{1.2}
\begin{array}{ll}
7)\ &
{{[[[x_{2},\,x_{3}],\,[x_{4},\,x_{6}]],\,
    [[x_{6},\,[x_{3},\,x_{4}]],\,[x_{6},\,
    [x_{4},\,x_{5}]]]]}=0};\\
&
    [[[x_{4},\,x_{6}],\,[x_{3},\,[x_{1},\,x_{2}]]],\,
    [[[x_{4},\,x_{5}],\,[x_{3},\,[x_{1},\,x_{2}]]],\\
&
   \hskip 10.3em [[x_{4},\,x_{6}],\,[x_{3},\,[x_{1},\,x_{2}]]]]]=0.
\end{array}
$

$
\arraycolsep=0pt
\renewcommand{\arraystretch}{1.2}
\begin{array}{ll}
5)\ &
{{[[x_{1},\,x_{2}],\,[x_{2},\,x_{3}]]}\,{=}\,0};\\
&
  {{[[[x_{2},\,x_{3}],\,[x_{4},\,x_{6}]],\,
    [[x_{6},\,[x_{3},\,x_{4}]],\,
    [x_{6},\,[x_{4},\,x_{5}]]]]}\,{=}\,0};\\
&
{{[[[x_{2},\,x_{3}],\,[x_{4},\,x_{6}]],\,
    [[[x_{2},\,x_{3}],\,[x_{4},\,x_{5}]],\,
    [[x_{2},\,x_{3}],\,[x_{4},\,x_{6}]]]]}=0}.
\end{array}
$

$
\arraycolsep=0pt
\renewcommand{\arraystretch}{1.2}
\begin{array}{ll}
6)\ &
  {{[x_{4},\,[x_{3},\,x_{5}]]}-{[x_{5},\,[x_{3},\,x_{4}]]}=0};\quad
  {{[[x_{2},\,x_{3}],\,[x_{3},\,x_{4}]]}=0};\\
&
  {{[[x_{2},\,x_{3}],\,[x_{3},\,x_{5}]]}=0};\quad
  {{[[x_{3},\,x_{4}],\,[x_{4},\,x_{6}]]}=0};\quad
  {{[[x_{4},\,x_{5}],\,[x_{4},\,x_{6}]]}=0};\\
&
  {{[[[x_{4},\,x_{6}],\,[x_{3},\,[x_{1},\,x_{2}]]],\,
    [[[x_{1},\,x_{2}],\,[x_{3},\,x_{4}]],\,
    [[x_{1},\,x_{2}],\,[x_{3},\,x_{5}]]]]}=0}.
\end{array}
$

$ \arraycolsep=0pt
\renewcommand{\arraystretch}{1.2}
\begin{array}{ll}
3)\ &
   {{[[x_{1},\,x_{2}],\,[x_{2},\,x_{3}]]}=0};\quad
  {{[[x_{2},\,x_{3}],\,[x_{3},\,x_{4}]]}=0};\quad
  {{[[x_{3},\,x_{4}],\,[x_{4},\,x_{6}]]}=0};\\&
  {{[[x_{4},\,x_{5}],\,[x_{4},\,x_{6}]]}=0};\quad
  {{[[x_{3},\,x_{4}],\,[[x_{3},\,x_{4}],\,
  [x_{4},\,x_{5}]]]}=0};\\&
  {{[[x_{6},\,[x_{3},\,x_{4}]],\,
      [[x_{3},\,x_{4}],\,[x_{5},\,[x_{3},\,x_{4}]
      ]]]}=0.}
      \end{array}
      $

      $
       \arraycolsep=0pt
\renewcommand{\arraystretch}{1.2}
\begin{array}{ll}
4)\ &
  {{[x_{4},\,[x_{3},\,x_{5}]]}-{[x_{5},\,[x_{3},\,x_{4}]]}=0};\quad
  {{[[x_{1},\,x_{2}],\,[x_{2},\,x_{3}]]}=0};\\
& {{[[x_{2},\,x_{3}],\,[x_{3},\,x_{4}]]}=0};\quad
  {{[[x_{2},\,x_{3}],\,[x_{3},\,x_{5}]]}=0};\\
& {{[[x_{3},\,x_{4}],\,[x_{4},\,x_{6}]]}=0};\quad
  {{[[x_{4},\,x_{5}],\,[x_{4},\,x_{6}]]}=0};\\&
  {{[[[x_{2},\,x_{3}],\,[x_{4},\,x_{6}]],\,
      [[x_{4},\,[x_{2},\,x_{3}]],\,[x_{5},\,
      [x_{2},\,x_{3}]]]]}=0.}
     \end{array}
     $

      $
       \arraycolsep=0pt
\renewcommand{\arraystretch}{1.2}
\begin{array}{ll}
8)\ &
  {{[[x_{2},\,x_{3}],\,[x_{3},\,x_{4}]]}=0};\quad
  {{[[x_{3},\,x_{4}],\,[x_{4},\,[x_{4},\,x_{5}]]]}=0};\\
& {{[[x_{4},\,x_{6}],\,[x_{4},\,[x_{4},\,x_{5}]]]}=0;}\\&
  {{[[x_{4},\,x_{5}],\,[x_{6},\,[x_{3},\,x_{4}]]]}-
  {2\,[[x_{4},\,x_{6}],\,[x_{5},\,[x_{3},\,x_{4}]]]}=0}.
  \end{array}
  $

$ \arraycolsep=0pt
\renewcommand{\arraystretch}{1.2}
\begin{array}{ll}
9)\ &
  {{[x_{4},\,[x_{3},\,x_{5}]]}-{[x_{5},\,[x_{3},\,x_{4}]]}=0};\\
& {{[[x_{1},\,x_{2}],\,[x_{2},\,x_{3}]]}=0};\quad
  {{[[x_{3},\,x_{4}],\,[x_{4},\,x_{6}]]}=0};\\&
  {{[[x_{4},\,x_{5}],\,[x_{4},\,x_{6}]]}=0};\quad
  {{[[x_{6},\,[x_{3},\,x_{4}]],\,
  [[x_{3},\,x_{4}],\,[x_{3},\,x_{5}]]
      ]}=0.}
      \end{array}
      $

      $
      \arraycolsep=0pt
\renewcommand{\arraystretch}{1.2}
\begin{array}{ll}
10)\ &
  {{[[x_{1},\,x_{2}],\,[x_{2},\,x_{3}]]}=0};\quad
  {{[[x_{2},\,x_{3}],\,[x_{3},\,x_{4}]]}=0};\\&
  {{[[x_{6},\,[x_{3},\,x_{4}]],\,
      [[x_{5},\,[x_{3},\,x_{4}]],\,
      [x_{6},\,[x_{3},\,x_{4}]]]]}=0};\\&
      {{[[x_{6},\,[x_{3},\,x_{4}]],\,
      [[x_{6},\,[x_{3},\,x_{4}]],\,
      [x_{6},\,[x_{4},\,x_{5}]]]]}=0.}
     \end{array}
     $

     $
       \arraycolsep=0pt
\renewcommand{\arraystretch}{1.2}
\begin{array}{ll}
11)\ &
  {{[[x_{2},\,x_{3}],\,[x_{3},\,x_{4}]]}=0};\quad
  {{[[x_{2},\,x_{3}],\,[x_{3},\,x_{5}]]}=0};\\&
  {{[[x_{6},\,[x_{3},\,x_{4}]],\,
  [[x_{3},\,x_{4}],\,[x_{3},\,x_{5}]]
      ]}=0};\\&
      {{[[[x_{2},\,x_{3}],\,[x_{4},\,x_{6}]],\,
      [[x_{4},\,[x_{2},\,x_{3}]],\,
      [x_{5},\,[x_{2},\,x_{3}]]]]}
     =0.}
     \end{array}
     $

     $
      \arraycolsep=0pt
\renewcommand{\arraystretch}{1.2}
\begin{array}{ll}
12)\ &
  {{[[x_{1},\,x_{2}],\,[x_{2},\,x_{3}]]}=0};\quad
  {{[[x_{5},\,[x_{2},\,x_{3}]],\,
  [[x_{3},\,x_{5}],\,[x_{4},\,x_{5}]]
      ]}=0};\\&
      {{[[x_{6},\,[x_{4},\,x_{5}]],\,
      [[x_{3},\,x_{5}],\,[x_{4},\,x_{5}]]]}=0};\\&
  {{[[[x_{4},\,x_{6}],\,[x_{5},\,[x_{2},\,x_{3}]]],\,
      [[x_{5},\,[x_{2},\,x_{3}]],\,
      [[x_{2},\,x_{3}],\,[x_{4},\,x_{5}]]
       ]]}=0.}
       \end{array}
       $

       $
        \arraycolsep=0pt
\renewcommand{\arraystretch}{1.2}
\begin{array}{ll}
13)\ &
  {{[[x_{1},\,x_{2}],\,[x_{2},\,x_{3}]]}=0};\quad
  {{[[x_{2},\,x_{3}],\,[x_{3},\,x_{4}]]}=0};\\
& {{[[x_{2},\,x_{3}],\,[x_{3},\,x_{5}]]}=0;}\quad
  {{[[x_{6},\,[x_{3},\,x_{4}]],\,
    [[x_{3},\,x_{4}],\,[x_{3},\,x_{5}]]]}=0;}\\
& {{[[[x_{4},\,x_{6}],\,[x_{3},\,[x_{1},\,x_{2}]]],\,
    [[[x_{1},\,x_{2}],\,[x_{3},\,x_{4}]],\,
    [[x_{1},\,x_{2}],\,[x_{3},\,x_{5}]]]]}=0.}
\end{array}
       $

       $
       \arraycolsep=0pt
\renewcommand{\arraystretch}{1.2}
\begin{array}{ll}
14)\ &
  {{[[x_{5},\,[x_{2},\,x_{3}]],\,[[x_{3},\,x_{5}],\,[x_{4},\,x_{5}]]]}=0};\\&
  {{[[x_{6},\,[x_{4},\,x_{5}]],\,[[x_{3},\,x_{5}],\,[x_{4},\,x_{5}]]]}=0};\\&
    [[[x_{3},\,[x_{1},\,x_{2}]],\,[x_{6},\,[x_{4},\,x_{5}]]],\,[[[x_{1},\,x_{2}],\,[x_{3},\,x_{5}]],\,\\
& \hskip 10.3em [[x_{4},\,x_{5}],\,[x_{3},\,[x_{1},\,x_{2}]]]]]= 0.
\end{array}
       $

       $
        \arraycolsep=0pt
\renewcommand{\arraystretch}{1.2}
\begin{array}{ll}
15)\ &
  {{[[x_{2},\,x_{3}],\,[x_{3},\,x_{4}]]}=0};\quad
  {{[[x_{3},\,x_{4}],\,[x_{4},\,x_{5}]]}=0};\\
& {{[[x_{3},\,x_{4}],\,[x_{4},\,x_{6}]]}=0};\quad
  {{[[x_{4},\,x_{6}],\,[[x_{4},\,x_{5}],\,[x_{4},\,x_{6}]]]}=0};\\
& {{[[x_{6},\,[x_{3},\,x_{4}]],\,[[x_{4},\,x_{6}],\,
    [x_{6},\,[x_{4},\,x_{5}]]]]}=0.}
\end{array}
      $

      $
        \arraycolsep=0pt
\renewcommand{\arraystretch}{1.2}
\begin{array}{ll}
16)\ &
  {{[[x_{1},\,x_{2}],\,[x_{2},\,x_{3}]]}=0};\quad
  {{[[x_{2},\,x_{3}],\,[x_{3},\,x_{5}]]}=0};\\
&
{{[[x_{3},\,x_{5}],\,[[x_{3},\,x_{5}],\,[x_{4},\,x_{5}]]]}=0};\;\;
  {{[[x_{6},\,[x_{4},\,x_{5}]],\,[[x_{3},\,x_{5}],\,[x_{4},\,x_{5}]]]}=0;}\\&
  {{[[[x_{3},\,x_{5}],\,[x_{4},\,x_{6}]],\,[[x_{3},\,x_{5}],\,[x_{4},\,[x_{3},\,x_{5}]]]]}=0.}
\end{array}
      $

      $
      \arraycolsep=0pt
\renewcommand{\arraystretch}{1.2}
\begin{array}{ll}
17)\ &
  {{[[x_{1},\,x_{2}],\,[x_{2},\,x_{3}]]}=0};\\
& {{[[[x_{2},\,x_{3}],\,[x_{4},\,x_{6}]],\,[[x_{4},\,[x_{2},\,x_{3}]],\,
    [x_{5},\,[x_{2},\,x_{3}]]]]}=0};\\&
  {{[[[x_{4},\,x_{6}],\,[x_{3},\,[x_{1},\,x_{2}]]],\,
    [[[x_{1},\,x_{2}],\,[x_{3},\,x_{4}]],\,[[x_{1},\,x_{2}],\,[x_{3},\,x_{5}]]]]}=0.}
\end{array}
$

 $
      \arraycolsep=0pt
\renewcommand{\arraystretch}{1.2}
\begin{array}{ll}
18)\ &
  {{[x_{5},\,[x_{4},\,x_{6}]]}-{[x_{6},\,[x_{4},\,x_{5}]]}=0};\quad
  {{[[x_{3},\,x_{4}],\,[x_{4},\,x_{5}]]}=0};\\&
  {{[[x_{3},\,x_{4}],\,[x_{4},\,x_{6}]]}=0};\quad
  {{[[x_{6},\,[x_{3},\,x_{4}]],\,
  [[x_{4},\,x_{6}],\,[x_{5},\,x_{6}]]
      ]}=0.}
\end{array}
$

$
      \arraycolsep=0pt
\renewcommand{\arraystretch}{1.2}
\begin{array}{ll}
19)\ &
  {{[[x_{2},\,x_{3}],\,[x_{3},\,x_{5}]]}=0};\quad
  {{[[x_{3},\,x_{5}],\,[x_{5},\,[x_{4},\,x_{5}]]]}=0};\\&
  {{[[x_{5},\,[x_{4},\,x_{5}]],\,[x_{6},\,[x_{4},\,x_{5}]]]}=0};\\
& {{[[[x_{3},\,x_{5}],\,[x_{4},\,x_{5}]],\,[[x_{3},\,x_{5}],\,[x_{4},\,x_{6}]]]}=0.}
\end{array}
$

$
   \arraycolsep=0pt
\renewcommand{\arraystretch}{1.2}
\begin{array}{ll}
20)\ &
  {{[[x_{5},\,x_{6}],\,[[x_{4},\,x_{5}],\,
  [x_{5},\,x_{6}]]]}=0};\;\;
  {{[[x_{5},\,[x_{3},\,x_{4}]],\,
  [[x_{4},\,x_{5}],\,[x_{5},\,x_{6}]]
      ]}=0};\\&
      {{[[[x_{3},\,x_{4}],\,[x_{5},\,x_{6}]],\,
      [[x_{5},\,x_{6}],\,[x_{6},\,[x_{4},\,x_{5}]]]]}=0.}
      \end{array}
      $

      $
      \arraycolsep=0pt
\renewcommand{\arraystretch}{1.2}
\begin{array}{ll}
21)\ &
  {{[[x_{5},\,x_{6}],\,[x_{5},\,[x_{4},\,x_{5}]]]}=0};\quad
  {{[[x_{5},\,[x_{3},\,x_{4}]],\,
  [x_{5},\,[x_{4},\,x_{5}]]]}=0};\\&
  {{[[[x_{3},\,x_{4}],\,[x_{5},\,x_{6}]],\,
      [[x_{4},\,x_{5}],\,[x_{5},\,x_{6}]]]}=0.}
      \end{array}
      $

\item[]

\hskip2.5ex \underline{$\fg(8, 6), \;\sdim =(133|56)$}

\item[]

$
\arraycolsep=0pt
\renewcommand{\arraystretch}{1.2}
\begin{array}{ll}
1)\ &
{{[[x_{5},\,x_{6}],\,[x_{6},\,x_{7}]]}=0};\\
&
{{[[[x_{5},\,x_{6}],\,[x_{4},\,[x_{1},\,x_{3}]]],\,
    [[[x_{2},\,x_{4}],\,[x_{5},\,x_{6}]],\,
    [[x_{3},\,x_{4}],\,[x_{5},\,x_{6}]]]]}=0}.
  \end{array}
$

$
\arraycolsep=0pt
\renewcommand{\arraystretch}{1.2}
\begin{array}{ll}
2)\ &
{{[[x_{4},\,x_{5}],\,[x_{5},\,x_{6}]]}=0};\quad
  {{[[x_{5},\,x_{6}],\,[x_{6},\,x_{7}]]}=0};\\
&
{{[[[x_{1},\,x_{3}],\,[x_{4},\,x_{5}]],\,
    [[x_{5},\,[x_{2},\,x_{4}]],\,[x_{5},\,
    [x_{3},\,x_{4}]]]]}=0}.
\end{array}
$

$
\arraycolsep=0pt
\renewcommand{\arraystretch}{1.2}
\begin{array}{ll}
3)\ &
    [[[x_{4},\,[x_{1},\,x_{3}]],\,
    [x_{7},\,[x_{5},\,x_{6}]]],\,[[[x_{6},\,x_{7}],\,[x_{5},\,[x_{2},\,x_{4}]]],\\
  & \hskip 12.5em [[x_{6},\,x_{7}],\,[x_{5},\,[x_{3},\,x_{4}]]]]]=0.
\end{array}
$

$ \arraycolsep=0pt
\renewcommand{\arraystretch}{1.2}
\begin{array}{ll}
4)\ &
  {{[[x_{2},\,x_{4}],\,[x_{4},\,x_{5}]]}=0};\quad
  {{[[x_{3},\,x_{4}],\,[x_{4},\,x_{5}]]}=0};\quad
  {{[[x_{4},\,x_{5}],\,[x_{5},\,x_{6}]]}=0};\\&
  {{[[x_{4},\,[x_{1},\,x_{3}]],\,
  [[x_{2},\,x_{4}],\,[x_{3},\,x_{4}]]
      ]}=0.}
      \end{array}
      $

      $
      \arraycolsep=0pt
\renewcommand{\arraystretch}{1.2}
\begin{array}{ll}
5)\ &
  {{[x_{3},\,[x_{2},\,x_{4}]]}-{[x_{4},\,[x_{2},\,x_{3}]]}=0};\quad
  {{[[x_{1},\,x_{3}],\,[x_{2},\,x_{3}]]}=0};\\
& {{[[x_{1},\,x_{3}],\,[x_{3},\,x_{4}]]}=0};\quad
  {{[[x_{2},\,x_{4}],\,[x_{4},\,x_{5}]]}=0};\quad
  {{[[x_{3},\,x_{4}],\,[x_{4},\,x_{5}]]}=0.}
  \end{array}
  $

  $
    \arraycolsep=0pt
\renewcommand{\arraystretch}{1.2}
\begin{array}{ll}
6)\ &
  {{[[x_{2},\,[x_{1},\,x_{3}]],\,[[x_{2},\,x_{3}],\,[x_{2},\,x_{4}]]]}=0};\\
& {{[[x_{5},\,[x_{2},\,x_{4}]],\,[[x_{2},\,x_{3}],\,[x_{2},\,x_{4}]]]}=0.}
\end{array}
      $

      $
      \arraycolsep=0pt
\renewcommand{\arraystretch}{1.2}
\begin{array}{ll}
7)\ &
  {{[[x_{1},\,x_{3}],\,[x_{2},\,x_{3}]]}=0};\quad
  {{[[x_{1},\,x_{3}],\,[x_{3},\,x_{4}]]}=0};\\&
  {{[[x_{5},\,[x_{3},\,x_{4}]],\,
  [[x_{2},\,x_{3}],\,[x_{3},\,x_{4}]]
      ]}=0.}
      \end{array}
      $

\enlargethispage{14mm}

      $
       \arraycolsep=0pt
\renewcommand{\arraystretch}{1.2}
\begin{array}{ll}
8)\ &
  {{[[[x_{1},\,x_{3}],\,[x_{4},\,x_{5}]],\,
      [[x_{2},\,[x_{1},\,x_{3}]],\,
      [x_{4},\,[x_{1},\,x_{3}]]]]}=0.}
     \end{array}
\vspace{-3mm}
     $

\end{arab}

}


{\appendix

\begin{landscape}
\vspace*{-1.5\baselineskip}
\section{The maximal root and the coefficients of linear dependence over $\Zee$ with
respect to the simple roots}

$
\arraycolsep=0pt
\renewcommand{\arraystretch}{1.1}
\begin{array}{ll}
 &\underline{\fg(2, 3), \;\sdim =(11|14)}\\
 1)\ & x_{11}=[[x_{2},\,x_{3}],\,[[x_{1},\,x_{2}],\,
 [x_{1},\,x_{3}]]];\quad
 \mathrm{weight}(x_{11})=(0,0,2); \quad(2, 2, 2).\\
  2)\ & x_{11}=
   [[x_{2},\,[x_{1},\,x_{3}]],\,[[x_1,\,x_{3}],\,[x_{2},\,x_{3}]]];\quad
   \text{weight} (x_{11})=(0 ,0 ,2); \quad(2, 2, 3).\\
 3)\ & x_{11}=
   [[x_{2},\,[x_{1},\,x_{3}]],\,[x_{3},\,[x_{2},\,x_{3}]]];\quad
   \text{weight} (x_{11})=(0 ,0 ,1 ); \quad(1, 2, 3).\\
 4)\ &
 x_{11}=[[x_{2},\,[x_{1},\,x_{3}]],\,[x_{3},\,[x_{1},\,x_{3}]]];\quad
        \text{weight} (x_{11})=(0 ,0 ,1 ); \quad(2, 1, 3).\\
 5)\ &
  x_{11}=[[x_{3},\,x_{3}],\,[x_{2},\,[x_{1},\,x_{3}]]];\quad
         \text{weight} (x_{11})=(0 ,0 ,1 ); \quad(1, 1, 3).\\[2mm]
 &\underline{\fg(1,6), \;\sdim =(21|14)}\\
 1)\ &
 x_{16}= [[x_{3},\,[x_{2},\,x_{2}]],\,
         [[x_{1},\,x_{2}],\,[x_{3},\,[x_{1},\,x_{2}]]]];\quad
         \text{weight}(x_{16})=(0,0,1); \quad(2, 4, 2).\\
\end{array}
$
\end{landscape}

\begin{landscape}
\vspace*{-1.5\baselineskip}
\[
\arraycolsep=0pt
\renewcommand{\arraystretch}{1.1}
\begin{array}{ll}
 &\underline{\fg(3, 3), \;\sdim =(21|16)}\\
 1)\ & x_{17}=
    [[x_{4},\,[x_{2},\,x_{3}]],\,
    [x_{3},\,[x_{3},\,[x_{1},\,x_{2}]]]];\quad
    \text{weight}(x_{17})=( 0, 0,1 ,0 ); \quad (1,2,3,1).\\
 2)\ & x_{17}=[[[x_{1},\,x_{2}],\,[x_{3},\,x_{4}]],\,
    [[x_{2},\,x_{3}],\,[x_{3},\,x_{4}]]];\quad
    \text{weight}(x_{17})=(0,0,1,0); \quad (1,2,3,2).\\
 3)\ & x_{17}= [[x_{4},\,[x_{2},\,x_{3}]],\,
    [[x_{1},\,x_{2}],\,[x_{3},\,x_{4}]]];\quad
    \text{weight} (x_{17})=( 0,0 ,2 ,0 ); \quad (1,2,2,2).\\
 4)\ & x_{17}=[[x_{4},\,[x_{3},\,[x_{1},\,x_{2}]]],\,
    [[x_{2},\,x_{3}],\,[x_{2},\,x_{4}]]];\quad
    \text{weight}(x_{17})=( 0,0 ,1 ,0 ); \quad (1,3,2,2).\\
 5) \ & x_{17}=[[x_{4},\,[x_{3},\,[x_{1},\,x_{2}]]],\,
  [[x_{2},\,x_{4}],\,[x_{3},\,[x_{1},\,x_{2}]]]];\quad
  \text{weight}(x_{17})=( 0,0 ,1 ,0 ); \quad (2,3,2,2).\\
 6) \ & x_{17}=
 [[x_{4},\,[x_{2},\,x_{3}]],\,
 [[x_{1},\,x_{2}],\,[x_{2},\,x_{3}]]
  ];\quad
  \text{weight}(x_{17})=(0,0,1,0); \quad (1,3,2,1).\\
7)\ & x_{17}=[[x_{4},\,[x_{3},\,[x_{1},\,x_{2}]]],\,
  [[x_{1},\,x_{2}],\,[x_{2},\,x_{3}]]];\quad
  \text{weight}(x_{17})=(0,0,1,0); \quad (2,3,2,1).\\
\end{array}
\]
\end{landscape}

\begin{landscape}
\vspace*{-2.5\baselineskip}
\[
\arraycolsep=0pt
\renewcommand{\arraystretch}{1.1}
\begin{array}{ll}
  & \underline{\fg(2, 6), \; \sdim =(35|20)}\\
 1)\ & x_{25}=[[[x_{2},\,x_{3}],\,[x_{3},\,x_{4}]],\,
    [[x_{1},\,x_{2}],\,[x_{5},\,[x_{3},\,x_{4}]]]];\quad
    \text{weight}(x_{25})=(0,0,0,1,0); \quad (1,2,3,2,1).\\
 2)\ & x_{25}=[[[x_{1},\,x_{2}],\,[x_{3},\,x_{4}]],\,
  [[x_{2},\,x_{4}],\,[x_{3},\,x_{5}]]];\quad
    \text{weight}(x_{25})=(0,0,0,-1,0); \quad (1,2,2,2,1).\\
 3)\ &x_{25}= [[[x_{1},\,x_{2}],\,[x_{5},\,[x_{3},\,x_{4}]]],\,
  [[x_{3},\,x_{5}],\,[x_{4},\,[x_{2},\,x_{3}]]]];  \quad
    \text{weight}(x_{25})=(0,0,0,-1,0); \quad (1,2,3,2,2).\\
 4)\ &x_{25}=[[[x_{2},\,x_{4}],\,[x_{3},\,x_{5}]],\,
  [[x_{2},\,x_{4}],\,[x_{3},\,[x_{1},\,x_{2}]]]];\quad
    \text{weight}(x_{25})=(0,0,0,1,0); \quad (1,3,2,2,1).\\
 5)\ & x_{25}=[[[x_{2},\,x_{4}],\,
 [x_{3},\,x_{5}]],\,
  [[x_{3},\,x_{4}],\,[x_{4},\,
  [x_{1},\,x_{2}]]]];\quad
    \text{weight}(x_{25})=(0,0,0,-1,0); \quad (1,2,2,3,1).\\
 6)\ & x_{25}=[[[x_{2},\,x_{4}],\,
 [x_{3},\,[x_{1},\,x_{2}]]],\,
  [[x_{3},\,x_{5}],\,[x_{4},\,
  [x_{1},\,x_{2}]]]];\quad
    \text{weight}(x_{25})=(0,0,0,1,0); \quad (2,3,2,2,1).\\[1mm]
 & \underline{\fg(8, 3), \;\sdim=(55|50)}\\
 1)\ & x_{50}=[[[[x_{1},\,x_{2}],\,[x_{3},\,x_{4}]],\,
    [[x_{2},\,x_{3}],\,[x_{4},\,x_{5}]]],\,
    [[[x_{3},\,x_{4}],\,[x_{3},\,[x_{1},\,x_{2}]]],\,
    [[x_{4},\,x_{5}],\,[x_{3},\,[x_{2},\,x_{3}]]]]];\\
    &\text{weight}(x_{50})=(0,0,0,0,1); \quad (2,4,6,4,2).\\
 2)\ & x_{50}=[[[[x_{2},\,x_{3}],\,[x_{4},\,x_{5}]],\,
   [[x_{4},\,x_{5}],\,[x_{3},\,[x_{1},\,x_{2}]]]],
   [[[x_{3},\,x_{4}],\,[x_{3},\,[x_{1},\,x_{2}]]],\,
   [[x_{4},\,x_{5}],\,[x_{3},\,[x_{2},\,x_{3}]]]]];\\
   & \text{weight}(x_{50})=(0,0,0,0,2); \quad (2,4,6,4,3).\\
3) \ &x_{50}=[[[[x_{2},\,x_{3}],\, [x_{5},\,[x_{3},\,x_{4}]]],\,
   [[x_{4},\,x_{5}],\,
   [x_{3},\,[x_{1},\,x_{2}]]]],\,
  [[[x_{3},\,x_{4}],\,
  [x_{4},\,[x_{2},\,x_{3}]]],\,
   [[x_{4},\,x_{5}],\,
   [x_{3},\,[x_{1},\,x_{2}]]]]];\\
 &\text{weight}(x_{50})=(0,0,0,0,1); \quad (2,4,6,5,5,3). \\
 4)\ & x_{50}=[[[[x_{1},\,x_{2}],\,[x_{4},\,x_{5}]],\,
   [[x_{2},\,x_{4}],\,[x_{3},\,x_{4}]]],\,
  [[[x_{2},\,x_{3}],\,[x_{4},\,x_{5}]],\,
   [[x_{4},\,x_{5}],\,
   [x_{3},\,[x_{1},\,x_{2}]]]]]\\
\end{array}
\]
    \end{landscape}
    \begin{landscape}
\vspace*{-2\baselineskip}
\[
\begin{array}{ll}
9)\ &
x_{50}=[[[[x_{1},\,x_{2}],\,[x_{4},\,x_{5}]],\,
   [[x_{2},\,x_{3}],\,[x_{2},\,x_{4}]]],\,
  [[[x_{2},\,x_{3}],\,[x_{4},\,x_{5}]],\,
   [[x_{2},\,x_{4}],\,[x_{3},\,[x_{1},\,x_{2}]]]]]\\
   &\text{weight}(x_{50})=(0,0,0,0,2); \quad (2,6,3,4,2). \\
10) \ & x_{50}=[[[[x_{2},\,x_{3}],\, [x_{4},\,x_{5}]],\,
   [[x_{1},\,x_{3}],\,
   [x_{4},\,[x_{1},\,x_{2}]]]],\,
  [[[x_{1},\,x_{2}],\,
  [x_{4},\,x_{5}]],\,
   [[x_{3},\,[x_{1},\,x_{2}]],\,
   [x_{5},\,[x_{2},\,x_{4}]]]]]\\
&\text{weight}(x_{50})=(0,0,0,0,2); \quad (4,5,3,4,3). \\ 11) \
&x_{50}= [[[[x_{2},\,x_{3}],\,[x_{4},\,x_{5}]],\,
   [[x_{1},\,x_{2}],\,[x_{3},\,[x_{1},\,x_{2}]]]],\,
  [[[x_{1},\,x_{2}],\,[x_{4},\,x_{5}]],\,
   [[x_{4},\,[x_{1},\,x_{2}]],\,
   [x_{4},\,[x_{2},\,x_{3}]]]]]\\
   &\text{weight}(x_{50})=(0,0,0,0,2); \quad (4,6,3,4,2). \\
   12) \ & x_{50}=
   [[[[x_{4},\,x_{5}],\,
   [x_{2},\,[x_{1},\,x_{3}]]],\,
   [[x_{4},\,x_{5}],\,
   [[x_{1},\,x_{3}],\,[x_{2},\,x_{3}]]]],\\
& \hskip 6.5ex [[[x_{4},\,x_{5}],\,
  [x_{2},\,[x_{1},\,x_{3}]]],\,
   [[x_{2},\,[x_{1},\,x_{3}]],\,
   [x_{4},\,[x_{2},\,x_{3}]]]]],\quad
   \text{weight}(x_{50})=(0,0,0,0,2); \quad (4,5,6,4,3). \\
\end{array}
    \]
    \end{landscape}
 \begin{landscape}
\vspace*{-1\baselineskip}
\[
\arraycolsep=0pt
\renewcommand{\arraystretch}{1.1}
\begin{array}{ll}
 &\underline{\fg(4, 6), \;\sdim =(66|32)}\\
 1)\ & x_{46}= [[[x_{6},\,[x_{3},\,x_{4}]],\,
   [[x_{2},\,x_{3}],\,[x_{4},\,x_{5}]]],\,
   [[x_{6},\,[x_{4},\,x_{5}]],\,[[x_{1},\,x_{2}],\,
   [x_{3},\,x_{4}]]]];\\
   & \text{weight}(x_{46})=(0,0,0,0,0,2); \quad (1,2,3,4,2,2).\\
 2) \ & x_{46}=[[[x_{5},\,[x_{2},\,x_{3}]],\,
 [x_{6},\,[x_{3},\,x_{4}]]],\,
  [[x_{6},\,[x_{4},\,x_{5}]],\,
  [[x_{1},\,x_{2}],\,[x_{3},\,x_{4}]]
   ]];\\
  & \text{weight}(x_{46})=(0,0,0,0,0,1);\quad (1,2,3,3,2,2).\\
3) \ & x_{46}=[[[x_{6},\,[x_{4},\,x_{5}]],\,
[[x_{2},\,x_{3}],\,[x_{4},\,x_{6}]]
   ],\,[[x_{6},\,[x_{5},\,[x_{3},\,x_{4}]]],\,
   [[x_{1},\,x_{2}],\,[x_{3},\,x_{4}]]]];\\
 & \text{weight}(x_{46})=(0,0,0,0,0,2);\quad (1,2,3,4,2,3).\\
 4) \ & x_{46}=
 [[[x_{6},\,[x_{3},\,x_{4}]],\,
 [x_{5},\,[x_{4},\,[x_{2},\,x_{3}]]]
   ],\,[[x_{6},\,[x_{3},\,x_{4}]],\,
   [[x_{1},\,x_{2}],\,[x_{3},\,x_{5}]]]];\\
   & \text{weight}(x_{46})=(0,0,0,0,0,1);\quad (1,2,4,3,2,2).\\
5) \ & x_{46}= [[[x_{6},\,[x_{4},\,x_{5}]],\,
[[x_{2},\,x_{3}],\,[x_{4},\,x_{5}]]
   ],\,[[[x_{1},\,x_{2}],\,[x_{3},\,x_{5}]],\,
   [[x_{3},\,x_{5}],\,[x_{4},\,x_{6}]]]];\\
   & \text{weight}(x_{46})=(0,0,0,0,0,1);\quad (1,2,3,3,4,2).\\
6) \ & x_{46}=[[[x_{6},\,[x_{3},\,x_{4}]],\,
[[x_{1},\,x_{2}],\,[x_{3},\,x_{5}]]
   ],\,[[x_{5},\,[x_{4},\,[x_{2},\,x_{3}]]],\,
   [[x_{2},\,x_{3}],\,[x_{4},\,x_{6}]]]];\\
   & \text{weight}(x_{46})=(0,0,0,0,0,1);\quad (1,3,4,3,2,2).\\
7) \ & x_{46}=[[[[x_{1},\,x_{2}],\,[x_{3},\,x_{4}]],\,
   [[x_{3},\,x_{5}],\,[x_{4},\,x_{6}]]],\,
  [[[x_{1},\,x_{2}],\,[x_{3},\,x_{5}]],\,
   [[x_{2},\,x_{3}],\,[x_{4},\,x_{6}]]]];\\
& \text{weight}(x_{46})=(0,0,0,0,0,1);\quad (2,3,4,3,2,2).\\
\end{array}%
\]
\end{landscape}

\section{The inverse matrices of the Cartan matrices}

In what follows the numbering of inverse matrices matches that of
the Cartan matrices. For the algebras not listed here, the Cartan
matrices are degenerate.

{\tiny
\[
\normalsize{\underline{\fg(1, 6)}} \tiny \quad 1)\left
(%
\]
}

}

\vfill


\footnotesize

\begin{thebibliography}{9999}

\bibitem[BjL]{BjL}
Bouarroudj S., Leites D., Simple Lie superalgebras and
non-integrable distributions in characteristic $p$. Zapiski
nauchnyh seminarov POMI, t. 331 (2006), 15--29 (F.A.Berezin
memorial volume); Reprinted in J. Math. Sci. (NY);
math.RT/0606682math.


\bibitem[CE]{CE}
Cunha I., Elduque A., An extended Freudenthal magic square in
characteristic $3$; math.RA/0605379

\bibitem[CE2]{CE2}
Cunha I., Elduque, A., The extended Freudenthal Magic Square and
Jordan algebras; math.RA/0608191


\bibitem[El1]{El1}
Elduque, A. New simple Lie superalgebras in characteristic 3. J.
Algebra 296 (2006), no. 1, 196--233

\bibitem[Er]{Er}
Ermolaev, Yu. B. Integral bases of classical Lie algebras.
(Russian) Izv. Vyssh. Uchebn. Zaved. Mat. 2004, , no. 3, 16--25;
translation in Russian Math. (Iz. VUZ) 48 (2004), no. 3, 13--22.

\bibitem[Eg]{Eg}
Egorov G. How to superize $\fgl(\infty)$. In: J.~Mickelsson e.a.,
(eds.) {\em Proc. Topological and Geometrical Methods in Field
Theory}, World Sci., Singapore, 1992, 135--146

\bibitem[FLS]{FLS}
Feigin B., Leites D., Serganova V., Kac--Moody superalgebras. In:
Markov M. et al (eds) {\it Group--theoretical methods in physics}
(Zvenigorod, 1982), v.  1, Nauka, Moscow, 1983, 274--278 (Harwood
Academic Publ., Chur, 1985, Vol.  1--3 , 631--637)

\bibitem[FSS]{FSS}
Frappat L., Sciarrino A.,
Sorba P., {\em Dictionary on Lie Superalgebras}. With 1 CD-ROM
(Windows, Macintosh and UNIX). Academic Press, Inc., San Diego,
CA, 2000. xxii+410 pp hep-th/9607161


\bibitem[FH]{FH}
Fulton, W., Harris, J., {\em Representation theory. A first
course}. Graduate Texts in Mathematics, 129. Readings in
Mathematics. Springer-Verlag, New York, 1991. xvi+551 pp


\bibitem[GG]{GG}
Grishkov, A.; Guerreiro, M. New simple Lie algebras over fields of
characteristic 2. Resenhas 6 (2004), no. 2-3, 215--221


\bibitem[Gr]{Gr} Grozman P., {\bf SuperLie},
\texttt{http://www.equaonline.com/math/SuperLie}


\bibitem[GL1]{GL1}
Grozman P., Leites D., Defining relations for classical Lie
superalgebras with Cartan matrix,  Czech.  J. Phys., Vol. 51,
2001, no.  1, 1--22; arXiv: \texttt{hep-th/9702073}


\bibitem[GL4]{GL4}
Grozman P., Leites D., Structures of $G(2)$ type and nonintegrable
distributions in characteristic $p$. Lett. Math. Phys.  74 (2005),
no. 3, 229--262; arXiv: \texttt{math.RT/0509400}

\bibitem[GLS]{GLS}
Grozman P., Leites D., Shchepochkina I., Invariant operators on
supermanifolds and standard models. In: In: M.~Olshanetsky,
A.~Vainstein (eds.)  {\em Multiple facets of quantization and
supersymmetry.  Michael Marinov Memorial Volume}, World Sci.
Publishing, River Edge, NJ, 2002, 508--555. [math.RT/0202193; ESI
preprint 1111 (2001)].

\bibitem[J1]{J1}
Jurman, G.,  A family of simple Lie algebras in characteristic
two.  J. Algebra 271 (2004), no. 2, 454--481.

\bibitem[K]{K}
Kac, V. {\em Infinite-dimensional Lie algebras}. Third edition.
Cambridge University Press, Cambridge, 1990. xxii+400 pp.


\bibitem[KKCh]{KKCh}
Kirillov, S. A.; Kuznetsov, M. I.; Chebochko, N. G. Deformations
of a Lie algebra of type $G\sb 2$ of characteristic three.
(Russian) Izv. Vyssh. Uchebn. Zaved. Mat. 2000, , no. 3, 33--38;
translation in Russian Math. (Iz. VUZ) 44 (2000), no. 3, 31--36

\bibitem[KL]{KL}
Kochetkov Yu., Leites D., Simple finite dimensional Lie algebras
in characteristic 2 related to superalgebras and on a notion of
finite simple group.  In: L.~A.~Bokut, Yu.~L.~Ershov and
A.~I.~Kostrikin (eds.)  {\em Proceedings of the International
Conference on Algebra.  Part 1., Novosibirsk, August 1989},
Contemporary Math. 131, Part 1, 1992, AMS, 59--67

\bibitem[KS]{KS}
Kostrikin, A. I., Shafarevich,  I.R., Graded Lie algebras of
finite characteristic, Izv. Akad. Nauk. SSSR Ser. Mat. 33 (1969)
251--322 (in Russian); transl.: Math. USSR Izv. 3 (1969) 237--304

\bibitem[KuCh]{KuCh}
Kuznetsov, M. I.; Chebochko, N. G. Deformations of classical Lie
algebras. (Russian) Mat. Sb. 191 (2000), no. 8, 69--88;
translation in Sb. Math. 191 (2000), no. 7-8, 1171--1190


\bibitem[Ku2]{Ku2}
Kuznetsov, M. I. Graded Lie algebras with the almost simple
component $L\sb 0$. Pontryagin Conference, 8, Algebra (Moscow,
1998). J. Math. Sci. (New York) 106 (2001), no. 4, 3187--3211.

\bibitem[Le1]{Le1}
Lebedev A., Non-degenerate bilinear forms in characteristic $2$,
related contact forms, simple Lie algebras and superalgebras.
arXiv: \texttt{math.AC/0601536}

\bibitem[Le2]{Le2}
Lebedev A., Simple Lie superalgebras preserving ortho-orthogonal
forms in characteristic $p=2$. IN PREPARATION


\bibitem[LL]{LL}
Lebedev A., Leites D., On realizations of the Steenrod algebras.
J. Prime Res. Math., v. 2, 2006,


\bibitem[LSh]{LSh}
Leites D., Shchepochkina I., Classification of the simple Lie
superalgebras of vector fields, preprint MPIM-2003-28
(\texttt{http://www.mpim-bonn.mpg.de})

\bibitem[vdL]{vdL}
Leur Johan van de., {\em Contragredient Lie superalgebras of
finite growth} (Ph.D. thesis) Utrecht, 1986; a short version
published in Commun.  in Alg., v.  17, 1989, 1815--1841

\bibitem[PS]{PS}
Penkov I., Serganova V., Generic irreducible representations of
finite dimensional Lie superalgebras.  Internat.  J. Math.  5,
1994, 389--419

\bibitem[Po]{Po}
Poletaeva E., The analogs of
Riemann and Penrose tensors on supermanifolds. math.RT/0510165

\bibitem[Pro]{Pro}
\text{ftp.mccme.ru/users/shuvalov/dyno/}


\bibitem[Sa]{Sa}
Sachse Ch., Sylvester-t'Hooft generators and relations between
them for $\fsl(n)$ and $\fgl(n|n$), Teor. Mat. Fiz 149(1), 2006,
3-17 (Russian; English translation in Theor. Math. Phys. 149(1),
2006, 1299-1311)

\bibitem[Se]{Se}
Serganova, V., Automorphisms of simple Lie superalgebras.
(Russian) Izv. Akad. Nauk SSSR Ser. Mat. 48 (1984), no. 3,
585--598

\bibitem[Se1]{Se1}
Serganova, V., On generalizations of root systems. Comm. Algebra
24 (1996), no. 13, 4281--4299



\bibitem[Shch]{Shch}
Shchepochkina I., How to realize Lie algebras by vector fields.
Theor. Mat. Fiz. 147 (2006) no. 3, 821--838; arXiv:
\texttt{math.RT/0509472}

\bibitem[Sh14]{Sh14}
Shchepochkina I., Five exceptional simple Lie superalgebras of
vector fields and their fourteen regradings. Representation Theory
(electronic journal of AMS), v. 3, 1999, 3 (1999), 373--415;
arXiv: \texttt{hep-th/9702121}

\bibitem[Sk]{Sk}
Skryabin, S. M. New series of simple Lie algebras of
characteristic $3$. (Russian) Mat. Sb. 183 (1992), no. 8, 3--22;
translation in Russian Acad. Sci. Sb. Math. 76 (1993), no. 2,
389--406

\bibitem[S]{S}
Strade, H. {\em Simple
Lie algebras over fields of positive characteristic. I. Structure
theory.} de Gruyter Expositions in Mathematics, 38. Walter de
Gruyter \& Co., Berlin, 2004. viii+540 pp.



\end{thebibliography}
\end{document}